\documentclass[a4paper,11pt]{article}
\usepackage{amsmath}
\usepackage{theorem}
\usepackage{amssymb,amsfonts,xypic}
\usepackage{a4wide}
\usepackage{color,ifthen}
\xyoption{all}
\newtheorem{theorem}{Theorem}[section]
\newtheorem{definition}[theorem]{Definition}
\newtheorem{proposition}[theorem]{Proposition}
\newtheorem{corollary}[theorem]{Corollary}
\newtheorem{lemma}[theorem]{Lemma}
\theorembodyfont{\rmfamily}
\newtheorem{note}{Note}
\newenvironment{proof}{\begin{trivlist}
\item[]{\em Proof\/}: }{\hspace*{\fill}
$\blacksquare$\end{trivlist}}

             \SelectTips{cm}{10}

\title{
             Postnikov Invariants
             of Crossed Complexes }
\author{
             M. Bullejos\footnote{Supported by DGI: BFM2001-2886.},
             E. Faro{$^*$},
             and
             M. A. Garc\'{\i}a-Mu\~{n}oz }

     \def\a{{\ensuremath{\mathcal A}}}
     
     \def\c{{\ensuremath{\mathcal C}}}
     
     \def\e{{\ensuremath{\mathcal E}}}
     \def\f{{\ensuremath{\mathcal F}}}
     \def\g{{\ensuremath{\mathcal G}}}

     \def\s{{\ensuremath{\mathcal S}}}
     \def\t{{\ensuremath{\mathcal T}}}

     \def\sff{{\ensuremath{\text{\sf F}}}}

     \def\bc{{\ensuremath{\mathbf C}}}
     
     \def\be{{\ensuremath{\mathbf E}}}
     \def\bf{{\ensuremath{\boldsymbol{f}}}}
     \def\bh{{\ensuremath{\boldsymbol{h}}}}
     
     \def\bg{{\ensuremath{\mathbf G}}}

     \newcommand{\cba}{{\ensuremath{\boldsymbol{\cal A}}}}
     
     \newcommand{\cbc}{{\ensuremath{\boldsymbol{\cal C}}}}
     
     \newcommand{\cbe}{{\ensuremath{\boldsymbol{\cal E}}}}
     \newcommand{\cbf}{{\ensuremath{\boldsymbol{\cal F}}}}
     \newcommand{\cbg}{{\ensuremath{\boldsymbol{\cal G}}}}

     \def\bbg{{\ensuremath{\mathbb G}}}

     \newcommand{\incl}{\text{\it i}}

     \newcommand{\Set}{{\ensuremath{\mathbf{Set}}}}
     \newcommand{\sets}{{\ensuremath{\mathbf{Set}}}}
     \newcommand{\Ab}{{\ensuremath{\mathbf{Ab}}}}
     \newcommand{\ab}{{\ensuremath{\mathbf{Ab}}}}
     \newcommand{\Top}{{\ensuremath{\mathbf{Top}}}}

     \newcommand{\xm}[1][]{{\ensuremath{\mathbf{Xm}_{#1}}}}

     \newcommand{\sgd}{{\ensuremath{\mathbf{SGd}}}}
     \newcommand{\acrs}{{\ensuremath{\mathbf{A}\crs{}}}}
     \newcommand{\scrs}[2][]{\ifthenelse{\equal{#1}{}}
                          {{\ensuremath{\mathbf{SCrs}_{#2}}}}
                          {{\ensuremath{\left(\mathbf{SCrs}_{#2}\right)_{#1}}}}}
     \newcommand{\sdel}{{\ensuremath{\boldsymbol{\Delta}}}}
     \newcommand{\asdel}{{\ensuremath{\bar{\boldsymbol{\Delta}}}}}
     \newcommand{\simpl}[2][]{\ifthenelse{\equal{#1}{}}
                          {{\ensuremath{\mathbf{Simp}\left(#2\right)}}}
                          {{\ensuremath{\mathbf{Simp}_{#1}\left(#2\right)}}}}
     \newcommand{\asimpl}[1]{\ifthenelse{\equal{#1}{}}
                         {{\ensuremath{\overline{\mathbf{Simp}}}}}
                         {{\ensuremath{\overline{\mathbf{Simp}}\left(#1\right)}}}}
     \newcommand{\dsimpl}[1]{\ifthenelse{\equal{#1}{}}
                         {{\ensuremath{\sets^{\op{\sdel}}}}}
                         {{\ensuremath{#1^{\op{\sdel}}}}}}
     \newcommand{\ddsimpl}[1]{\ifthenelse{\equal{#1}{}}
                         {{\ensuremath{\sets^{\op{\sdel}\times \op{\sdel}}}}}
                         {{\ensuremath{#1^{\op{\sdel}\times \op{\sdel}}}}}}
     \newcommand{\adsimpl}[1]{\ifthenelse{\equal{#1}{}}
                         {{\ensuremath{\sets^{\op{\asdel}}}}}
                         {{\ensuremath{#1^{\op{\asdel}}}}}}
     
     \newcommand{\agpd}{{\ensuremath{\mathbf{A}\Gpd}}}
     
     \newcommand{\pxm}[1][]{{\ensuremath{\mathbf{Pxm}_{#1}}}}

     \newcommand{\crs}{\ensuremath{\mathbf{Crs}}}
     \newcommand{\Graphs}{{\ensuremath{\mathbf{Gph}}}}
     \newcommand{\Gph}{\Graphs}

     \newcommand{\Groupoids}{{\ensuremath{\mathbf{Gpd}}}}
     \newcommand{\Gpd}{\Groupoids}
     
     \newcommand{\tdGpd}{{\ensuremath{\mathbf{TdGpd}}}}
     \newcommand{\Groups}{{\ensuremath{\mathbf{Gr}}}}
     \newcommand{\Gr}{\Groups}
     \newcommand{\Gp}{\Groups}
     \newcommand{\gp}{\Groups}

     \newcommand{\sset}{{\ensuremath{\boldsymbol{\mathcal{S}}}}}

     \newcommand{\op}[1]{{\ensuremath{#1^{\text{\rm\kern-0.1em o\kern-0.1em p}}}}}
     \newcommand{\cte}[1]{{\ensuremath{#1^{\scriptscriptstyle\text{ct}}}}}
     \newcommand{\initxm}[1][\vphantom{.}]{{\boldsymbol0}_{_{#1}}}
     \newcommand{\termxm}[1][\vphantom{.}]{{\boldsymbol1}_{_{#1}}}
     \def\cotrsr{{\cotrip[n]^\bullet}}
     \def\semi#1#2{{\ensuremath{#1\ltimes #2}}}
     \def\actor(#1){{\ensuremath{\text{Act}(#1)}}}
     \def\scoH{H_{\text{\rm sing}}}
     \def\actt(#1){{\ensuremath{\text{Act1}(#1)}}}
     \def\act2(#1){{\ensuremath{\text{Act2}(#1)}}}
     \def\aut(#1){{\ensuremath{\text{AUT}(#1)}}}
     \def\ner{{\ensuremath{\mathop{\rm Ner}}}}
     \def\naut(#1){{\ensuremath{\text{NAUT}(#1)}}}
     \def\Fib(#1){{\ensuremath{\text{\bfseries Fib}(#1)}}}
     \def\fib{{\ensuremath{\text{\sf fib}}}}

     \def\laction#1#2{{}^{#1}\!#2}
     
     \newcommand{\wbar}[1][]{{\ensuremath{\overline{W}_{\!#1}}}}
     \newcommand{\W}{{\ensuremath{\overline{W}}}}
     \newcommand{\cotrip}[1][]{\ensuremath{\bbg_{#1}}}
     \newcommand{\Img}{\mathop{\rm im}}
     
     \newcommand{\End}{\text{\sf End}}

     \newcommand{\Tor}{\text{\it Tor}}
     \newcommand{\Ext}{\text{\it Ext}}
     \newcommand{\fbr}{\text{\sf fbr}}
     
     \newcommand{\obj}{\text{\sf obj}}
     
     \newcommand{\arr}{\text{\sf arr}}
     \newcommand{\discr}{\text{\sf discr}}
     \newcommand{\codiscr}{\text{\sf codiscr}}
     \newcommand{\gpd}{\text{\sf gpd}}
     \newcommand{\fxm}{\text{\sf xm}}
     \newcommand{\base}{\text{\sf base}}
     
     \newcommand{\zero}{\text{\sf zero}}
     \newcommand{\ins}{\text{\sf ins}}
     \newcommand{\abgc}[2]{\tilde#1_{#2}}
     \newcommand{\abgp}[2]{\tilde#1_{#2}}
     
     \newcommand{\cosk}{\text{\sf cosk}}

     \renewcommand{\hbar}{{\ensuremath{\bar{h}}}}

     \def\myxymat#1{
                  \vcenter{
                  \xymatrix}}
     \def\adj(#1,#2,#3,#4){
     {\xymatrix{{#1:#2} \ar@<.6ex>[r] & {#3:#4} \ar@<.6ex>[l] }}}
     \def\xto#1{\xrightarrow{#1}}
     \def\lto{\leftarrow}
     \newcommand\hto{\hookrightarrow}

     \newcommand{\p}{Postnikov}
     
     \newcommand{\pto}{\p\ tower}
     
     \newcommand{\pin}{\p\ invariant}

     \newcommand{\cc}{crossed complex}
     \newcommand{\ccs}{\cc es}
     \newcommand{\ncc}[1][]{\ifthenelse{\equal{#1}{}}
                         {\ensuremath{n}}
                         {\ensuremath{#1}}-\cc}
     \newcommand{\nccs}[1][]{\ifthenelse{\equal{#1}{}}
                          {\ensuremath{n}}
                          {\ensuremath{#1}}-\ccs}
     \newcommand{\npcc}[1][]{\ifthenelse{\equal{#1}{}}
                       {\ensuremath{n}}
                       {\ensuremath{(#1)}}-\cc}
     \newcommand{\npccs}[1][]{\ifthenelse{\equal{#1}{}}
                         {\ensuremath{n}}
                         {\ensuremath{(#1)}}-\ccs}

     \def\quitar#1{}
     \def\guardar#1{}

     \def\leq{\leqslant}
     \def\geq{\geqslant}
     \def\tilde{\widetilde}

     \def\bpar#1{\ensuremath{\big(}#1\ensuremath{\big)}}

     \def\myxymat#1{\SelectTips{cm}{10}
                 \vcenter{\xymatrix@C=1.2pc@R=1.2pc{#1}}}

     \newcommand{\utor}{\underline{\Tor}}
     \def\h{{\ensuremath{\mathcal H}}}

     \def\techo{{\ensuremath{{\mathsf{ceil}}}}}
     \def\ntecho{{\ensuremath{{\techo_n}}}}

     \newcommand{\Ker}{\text{\sf ker}}

     \def\cuad(#1,#2,#3,#4)(#5,#6,#7,#8){                    
     \diagram                                                
     #1 \rto^{#5}\dto_{#8} & #2\dto^{#6} \\                  
     #4 \rto_{#7} & #3                                       
     \enddiagram                                             
     }
     \def\gcrs[#1]{{\ensuremath{\mathbf{G}\crs_{#1}}}}

     \def\N{{\ensuremath{\mathbf N}}}
     \def\fcrs{{\ensuremath{\text{\sf crs}}}}
     \def\ck{\ensuremath{\mathcal K}}

\begin{document}
\maketitle
\setcounter{tocdepth}{2}

\begin{abstract}
We determine the \p\ Tower and \p\ Invariants of a Crossed Complex in a purely algebraic way.
Using the fact that Crossed Complexes are homotopy types for filtered spaces, we use the above
``algebraically defined'' \pto\  and \pin{s} to obtain from them those of filtered spaces. We
argue that a similar ``purely algebraic'' approach to \pin{s} may also be used in other
categories of spaces.
\end{abstract}

\section{Introduction}

The theory of Postnikov towers provides both, a way of analyzing a space $X$
from the point of view of its homotopy groups, and a prescription for the
construction of spaces with specified homotopy groups in each dimension. The
required data for this construction is the information contained in the
Postnikov tower of the space: a diagram of spaces
$$
\xymatrix@C=1.5pc@R=1.35pc {\cdots
\ar[r]& X_{n+1}\ar[r]^-{\eta_{n+1}} & X_n\ar[r]^-{\eta_n} & X_{n-1}\ar[r]
&\cdots \ar[r]& X_0,}
$$
whose inverse limit has the same homotopy type as the given space and
where each map $\eta_n$ is a fibration whose fibers are Eilenberg-Mac
Lane spaces of the type of a $K(\Pi,n)$.

The Postnikov invariants of the space $X$ are cohomology invariants, denoted
$k_{n}$, $n\geq1$, which provide the necessary information in order to build the
Postnikov tower of $X$ floor by floor. The Postnikov invariant $k_{n}$,
associated to the fibration $\eta_n$, says how to glue $K(\Pi,n)$ spaces into
the space $X_{n-1}$ to form $X_{n}$.

The purpose of this paper is to present a purely algebraic approach to the
calculation of the Postnikov invariants of a space in the sense that it avoids
the use of topological tools such as universal covering; these are replaced by
algebraic tools such as free resolutions.

One of the motivations for such an algebraic approach is the fact that it offers
the possibility of applying it to more complicated contexts such as categories
of diagrams of spaces. Previous work in this direction can be seen in
\cite{BuCaFa1998}, where the third equivariant Postnikov invariant of a
$G$-space is calculated by purely algebraic methods of the same nature as those
presented here.

Our approach to Postnikov towers and Postnikov invariants of the spaces in a
given category of spaces \t\ is based in the existence of an algebraic category
$\s$ with a Quillen model structure, together with a pair of functors
$\Pi:\t\rightarrow\s$, $B:\s\rightarrow \t$ which  induce an
equivalence in the corresponding homotopy categories. Hence \s\ is a
category of algebraic models for the homotopy types of the spaces in
\t. In this situation, the calculation of the Postnikov towers of the
spaces in \t\ can be reduced to calculating Postnikov towers in \s\
provided that the ``classifying space'' functor
$B$ preserves fibrations as well as the homotopy type of their fibers.
This is the case for the functors which are the object of this paper.
These are, on the one hand, the functor $\Pi=${\it ``Fundamental
crossed complex of the singular complex of a space''}, and on the other
hand, the functor $B=${\it``Geometric realization of the nerve
of a crossed complex''} (see below).

The calculation of Postnikov towers of the algebraic models which are the
objects of \s\ is based in the following general scheme: For every non-negative
integer $n$ we seek a full, reflective subcategory $i_n:\s_n\rightarrow\s$ whose
objects model all ``\emph{homotopy $n$-types}'' in \s, and such that $\s_n$ is
contained in $\s_{n+1}$ in such a way that the inclusions $j_n:\s_n\to\s_{n+1}$
satisfy $i_{n+1}j_n=i_n$. Then if $\tilde{P}_n$ is the left adjoint to $i_n$,
the identity $\tilde{P}_{n+1}i_{n+1}=1_{\s_{n+1}}$ implies
$\tilde{P}_{n+1}i_n=j_n$ and the composites $P_n=i_n\tilde{P}_n$ are idempotent
endofunctors of \s\ verifying $P_{n+1}P_n=P_n\simeq P_nP_{n+1}$, and related by
a chain of natural transformations
\begin{equation}
\label{chain}
\xymatrix@C=1.5pc@R=1.35pc{\cdots
\ar[r]& P_{n+1}\ar[r]^-{\eta_{n+1}} & P_n\ar[r]^-{\eta_n} & P_{n-1}\ar[r]
&\cdots \ar[r]& P_0}
\end{equation} (where $\eta_{n+1}$ is the image by $P_{n+1}$ of the unit of the
adjunction $\tilde{P}_n\vdash i_n$). Then we prove that this chain is the
``\emph{universal Postnikov tower}" in \s, in the sense that for any object
$C\in\s$ the evaluation of the above chain \eqref{chain} in $C$ yields the
Postnikov tower of $C$.

With regards to the Postnikov invariants, it is noteworthy the simple form taken
by the fibrations of the Postnikov towers of the objects of \s, making it easy
to analyze them and to show that the component of $\eta_{n+1}$ in each object
can be interpreted as a 2-extension, a 2-torsor, and then it gives rise to an
element $k_{n+1}$ of a 2-dimensional algebraic (cotriple) cohomology in (a slice
of) the category $\s_n$ of algebraic $n$-types.  We regard such 2-dimensional
cohomology element as a sort of ``\emph{algebraic Postnikov invariant}", the
``\emph{topological}" one residing in a $(n+2)$-dimensional singular cohomology.

The last step in our approach consists in obtaining the topological Postnikov
invariants from the algebraic ones. This is achieved by showing the existence of
a natural map from the algebraic 2-cohomology of an algebraic $n$-type in $\s_n$
to the singular $(n+2)$-cohomology of its corresponding classifying space.

The ideal scenario offering the necessary tools to apply the algebraic approach
just described is that in which \t\ is the category of CW-complexes and \s\ is
the category of simplicial groupoids. The main interest of this context lies, of
course, in the fact that simplicial groupoids model all homotopy types and,
therefore, a procedure to calculate the algebraic Postnikov invariants of
simplicial groupoids could be used to obtain the Postnikov invariants of any
space. The work presented here is, however, more modest in scope and it is, in
fact, a preliminary step in that direction. We carry out the general method
described above in the category of crossed complexes, a category which does not
model all homotopy types. However, although our present results cannot be used
to obtain the Postnikov invariant of all spaces, they are, of course, sufficient
to obtain the Postnikov invariants of any space having the homotopy type of a
crossed complex.

The general plan of the paper is as follows: Section \ref{xm} serves to set-up
our notation and to introduce the definition and main facts about crossed
modules that are used in the paper. Everything here is review material which can
be found elsewhere in the literature, except that is presented in a, perhaps,
slightly non-conventional way, with an emphasis in the functorial aspect of the
definitions. We apologize for any distraction this may cause to those readers
who are already familiar with the subject. Section \ref{crs} introduces crossed
complexes, the categories of $n$-types in crossed complexes and the Postnikov
towers they give rise to. Crossed complexes are again introduced in a slightly
non conventional way, being defined in terms of crossed modules instead of in
terms of groupoids, as it is customary. Choosing a definition which is based on
a more elaborate concept not only simplifies the definition itself but, more
important, it allows simpler and clearer reasonings and proofs. We also show in
this section that the geometric realizations of these Postnikov towers are the
Postnikov towers of spaces. Section \ref{inv} is the main section of the paper.
Here the fibrations in the Postnikov towers of crossed complexes are analyzed
and interpreted as extensions, torsors, and therefore, as a consequence of
Duskin's interpretation theorem, as cohomology elements in a cotriple
cohomology. Finally a general theorem is proved showing how to map the cotriple
cohomology of crossed complexes to the singular cohomology of its classifying
spaces. Section \ref{apen} is an appendix containing the basic definitions and
results about torsors and their role in Duskin's interpretation theorem of
cotriple cohomology. This material, essential for the main results of the paper,
is well known to the specialist but it is no so well known in larger circles. It
has been put in an appendix in order not to break the discourse and to allow the
reader to focus on the main line of reasoning.

\section{Crossed modules}
\label{xm}

We denote \Gr\ the category of groups and \Gpd\ the category of small groupoids,
that is, the category of internal groupoids in the category $\sets$ of sets. By
\tdGpd\ we denote the full subcategory of $\Gpd$ determined by the totally
disconnected groupoids. If $X$ is a set, $\Gpd_{X}$ denotes the subcategory of
$\Gpd$ whose objects are all groupoids with set of objects $X$ and whose arrows
are functors which are the identity on objects. Similarly,  $\tdGpd_X$ denotes
the full subcategory of $\Gpd_X$ determined by the totally disconnected
groupoids. Clearly, $\tdGpd_X$ can be identified with the category $\Gr(\Set/X)$
of internal group objects in the slice category $\Set/X$. For a given groupoid
\bg\  we denote $\obj(\bg)$ its set of objects, and $\arr(\bg)$ its set of
arrows. It is clear that $\obj$ determines a functor $\obj:\tdGpd\rightarrow
\Set$ whose fiber over a set $X$ is the category $\tdGpd_X$.

If \bg\ is a groupoid, a (left) \bg-group is a functor from \bg\ to \Gr. We will
use exponential notation to denote functor categories, so that the category of
(always left-) \bg-groups will be denoted $\Gr^\bg$. An important example of a
\bg-group is the functor $\End_\bg:\bg\rightarrow\Gr$ taking each object of \bg\
to its group of endomorphisms and each arrow $u$ in \bg\ to the group
homomorphism (really an iso) given by conjugation by $u$. This \bg-group is
often referred to as the groupoid \bg\ acting on itself by conjugation.

A given \bg-group $C:\bg\rightarrow\Gr$ is often determined in terms of an
action of (the arrows of) \bg\ on (the arrows of) a totally disconnect groupoid,
$\widehat{C}$, whose set of object is $\obj(\bg)$ and whose endomorphism groups
are $\End_{\widehat{C}}(x)=C(x)$. This action of \bg\ on $\widehat{C}$ is
traditionally denoted
$$
\laction t u= C(t)(u),
$$ for $t:x\rightarrow y$ an arrow in \bg\ and $u$ an element in $C(x)$. This
description of \bg-groups is objectified by a full and faithful functor
$\widehat{(\;)}:\Gr^\bg\rightarrow\tdGpd_{\obj(\bg)}$ which reflects zero
objects and zero maps, and therefore not only preserves but also reflects chain
complexes. Obviously, $\widehat{\End_{\bg}} = \End(\bg)$, the subcategory of
$\bg$ consisting of just its endomorphisms.

For a given groupoid \bg\ we denote $\pxm_\bg$ the category of pre-crossed
modules over \bg, which we define as the slice category
$\pxm_\bg=\Gr^\bg/\End_\bg$. The initial and terminal objects in $\pxm_\bg$ are
denoted $\initxm[\bg]$ and $\termxm[\bg]$ respectively, so that $\initxm[\bg]$
is a constant zero functor $\bg\rightarrow\Gr$ together with the unique natural
transformation from it to $\End_\bg$, while $\termxm[\bg]$ is the functor
$\End_\bg$ together with its identity map. Note that $\initxm[\bg]=\termxm[\bg]$
if and only if \bg\ is discrete as category (that is, all arrows in \bg\ are
identities).

If \bg\ and $\bg'$ are any two groupoids and $(C,\delta)$,
$(C',\delta')$ are pre-crossed modules respectively over \bg\ and
$\bg'$, a morphism of pre-crossed modules from $(C,\delta)$ to
$(C',\delta')$ is a pair $(f,\alpha)$ where
$f:\bg\rightarrow\bg'$ is a \emph{change-of-base} functor and
$\alpha:C\rightarrow C'\circ f$ is a natural transformation such that
$(\delta'*f)\circ\alpha=\tilde{f}\circ\delta$, where $\tilde{f}$ is the
same functor $f$ but regarded as natural transformation from $\End_\bg$ to
$\End_{\bg'}\circ f$,
$$
\xymatrix@C=1.5pc@R=1.75pc {C\ar[d]_\delta
\ar[r]^-\alpha & C'\circ f\ar[d]^{\delta'*f} \\
\End_\bg\ar[r]_-{\tilde{f}} & \End_{\bg'}\circ f. }
$$ For an object $x\in\bg$ and an element $u\in C(x)$, this condition reads
$f(\delta_x(u))=\delta'_{f(x)}(\alpha_x(u))$. The general morphisms of
pre-crossed modules just defined are the arrows of the category of pre-crossed
modules, denoted \pxm. The structure of an object in \pxm\ can be described as a
triple $(\bg,C,\delta)$ where \bg\ is a groupoid and $(C,\delta)$ is a
pre-crossed module over \bg. By a \emph{reduced} pre-crossed module we mean one
in which \bg\ is just a group.

The following proposition provides the definition of the fundamental groupoid of a pre-crossed
module. Note that if $(C,\delta)$ is a pre-crossed module over \bg, by applying the functor
$\widehat{(\;)}$ to the natural map $\delta:C\rightarrow \End_\bg$ we obtain a functor
$\widehat{\delta}:\widehat{C}\rightarrow \End(\bg)$ which ($\widehat{C}$ being totally
disconnected) is equivalent to a functor $\widehat{C}\rightarrow\bg$. The latter will not be
distinguished from $\widehat{\delta}$.
\begin{proposition}
\label{adjs betw xm and gpd} The categories $\pxm_\bg$ are the fibres of a fibration ``base
groupoid of a pre-crossed module'', $\base:\pxm\rightarrow\Gpd$. This functor has both adjoints
$\discr\dashv\base\dashv\codiscr$ given by the initial (left) and terminal (right) objects in the
corresponding fibres. Furthermore the left adjoint $\discr$ has a further left adjoint
``\emph{fundamental groupoid}" $\pi_1\dashv \discr$.

\end{proposition}
\begin{proof} Everything is quite standard; we just comment on the last
statement. The fundamental groupoid of a pre-crossed module $\c=(\bg,C,\delta)$ is calculated by
the coequalizer
$$
\everyentry={\vphantom{\Big(}}
\xymatrix@C=1.75pc@R=1.75pc
{\widehat{C}\ar@<.5ex>[r]^{\widehat{\delta}}\ar@<-.5ex>[r]_0 &
\bg\ar[r]^-q & \pi_1(\c),}
$$ that is, the fundamental groupoid is given by the quotient
$\pi_1(\c)=\bg/\Img(\widehat{\delta})$. Note that all functors in the above
diagram are the identity on objects.
\end{proof}

For a given pre-crossed module $\c=(\bg,C,\delta)$ the functors
$C\circ\widehat\delta$ and $\End_{\widehat{C}}$ agree on objects but, in
general, not in arrows. Therefore, that these two functors be equal is a special
property a pre-crossed module may have.

\begin{definition} A crossed module is a pre-crossed module $\c=(\bg,C,\delta)$
such that $C\circ\widehat{\delta}=\End_{\widehat{C}}$. The category of crossed
modules, denoted \xm, is the corresponding full subcategory of \pxm. For a given
groupoid \bg, the category of \bg-crossed modules, denoted $\xm_\bg$, is the
obvious full subcategory of $\pxm_\bg$.
\end{definition}

In terms of elements, the condition that $C\circ\widehat{\delta}$ and
$\End_{\widehat{C}}$ agree on arrows reads:
$$
\laction {\delta_x(u)}v= u v u^{-1},
$$ for all objects $x\in\bg$ and elements $u,v\in C(x)$. This is the well known
Peiffer identity. This property implies that $\ker\delta_x$ is contained in the
center of $C(x)$, and this, in turn, has the following important consequences:
\begin{enumerate}
\item For any given \bg-group $C$, the \bg-pre-crossed module $(C,0)$ is a
crossed module if and only if every group $C(x)$ is abelian, that is, if $C$ is
a \bg-module.
\item For any \bg-crossed module $(C,\delta)$, the kernel of $\delta$
(calculated in $\Gr^\bg$) is a \bg-module (that is,
$(\ker\delta)(x)=\ker\delta_x$ is an abelian group).
\item The action of $\Img\delta$ on $\ker\delta$ is trivial, that is, the
following diagram commutes
\begin{equation}
\label{consequence xm}
\everyentry={\vphantom{\Big(}}
\xymatrix@C=1.75pc@R=1.75pc
{\widehat{C}\ar@<.6ex>[r]^{\widehat{\delta}}\ar@<-.4ex>[r]_0 & \bg
\ar[r]^-{\ker\delta} &
\Ab,}
\end{equation} where we have denoted ``0'' the functor which is the identity on
objects and sends every map to an identity; this functor, if regarded as a map
in $\Gpd_{\obj(\bg)}$, is indeed a zero map.
\end{enumerate}

\bigskip
\noindent {\bfseries Examples:}
\vskip0pt  1. For any $\bg$-module $A:\bg\to\Ab$, the pre-crossed module
$\zero(A) = (\bg,A,0)$ is a crossed module.
\vskip0pt 2. Any pre-crossed module $(\bg, C, \delta)$ with $\delta$ a
monomorphism is a crossed module.
\vskip0pt 3. In particular, for any groupoid $\bg$, the pre-crossed modules
$\initxm[\bg]$ and $\termxm[\bg]$ are crossed modules.

\bigskip As a consequence of the above Example 3, the right and left adjoints to
``base groupoid of a pre-crossed module'' are also right and left
adjoints to ``base groupoid of a crossed module'' and furthermore
``fundamental groupoid of a crossed module'' is left adjoint to
``discrete crossed module on a groupoid''. From now on we will
regard the functors in the sequence of adjunctions $\pi_1 \dashv
\discr \dashv \base \dashv \codiscr$ of Proposition \ref{adjs betw
xm and gpd} as defined/taking values in $\xm$.

There is an important forgetful functor defined on the category of \bg-crossed
modules, which will be used later on. This is the functor \begin{equation}
\label{techo}\techo_2=\techo:\xm_\bg\rightarrow \gp^\bg \end{equation} taking a
\bg-crossed module $\c=(\bg,C,\delta)$ to the \bg-group $C$ and each map
$(1_{\bg},\alpha):(\bg,C,\delta)\to(\bg,C',\delta')$ of \bg-crossed modules to
the natural transformation $\alpha:C\to C'$. It is an important property of this
functor the fact  that it preserves finite limits and coequalizers.

Our next objective is to establish the tripleability of the category of crossed
modules over a certain category (see below) so that we can define the cotriple
which will be used to calculate an algebraic cohomology of crossed modules.

\begin{proposition}
\label{free xm on a pxm} The inclusion functor $U:\xm\rightarrow\pxm$ has a left
adjoint which is calculated by factoring out the Peiffer subgroup.
\end{proposition}
\begin{proof} See \cite[p. 9]{BrWe1996}, \cite{BrHu1982}, or \cite{HoMeSi1993}.
\end{proof}

This inclusion functor $U:\xm\rightarrow\pxm$ is in fact monadic. We will use it
to obtain, by composing it with a certain forgetful functor
$U':\pxm\rightarrow\agpd$, another monadic functor which will determine the
cotriple on $\xm$ by means of which we will calculate the cohomology of crossed
modules. Let \agpd\ be the category of ``arrows to groupoids'' whose objects are
triples $(X,f,\bg)$ where $X$ is a set, \bg\ is a groupoid, and $f:X\rightarrow
\End(\bg)$ is a map from $X$ to the set of all arrows of \bg\ which are
endomorphisms. An arrow from $(X,f,\bg)$ to $(X',f',\bg')$ in
\agpd\ is a pair $(\alpha,\beta)$ where $\alpha:X\rightarrow X'$ is a map
of sets, and $\beta:\bg\rightarrow\bg'$ is a functor such that
$f'\alpha=\beta f$. Then we have:

\begin{proposition}
\label{pxm to agpd has ladj} The obvious forgetful functor
$U':\pxm\rightarrow\agpd$ has a left adjoint.
\end{proposition}
\begin{proof} The forgetful functor $U':\pxm\rightarrow\agpd$ takes a
pre-crossed module $(\bg,C,\delta)$ to the triple
$(\arr(\widehat{C}),\widehat{\delta},\bg)$. We will merely give the definition
of its left adjoint $F:\agpd\rightarrow\pxm$. This is defined on objects as
$F(X,f,\bg)=(\bg,C,\delta)$, where $C:\bg\rightarrow\Gr$ is defined on objects as
$$ C(x)= F_{\text{gp}} \bigg( \coprod_{z\in \obj(\bg)} \Big(  G(z,x)\times
\coprod_{v\in\bg(z,z)} \fbr(f,v) \Big) \bigg) ,
$$
$F_{\text{gp}}$ is the free group functor, and $\fbr(f,v)$ is the fiber of the
map $f$ at $v$. Thus, $C(x)$ is the free group generated by all pairs $\langle
t,u\rangle $, where $t:z\rightarrow x$ is a map in \bg\ and $u\in X$ is such
that $f(u)$ is an endomorphism of $z$ in \bg. The \bg-group $C$ is defined on
arrows $s:x\rightarrow y$ in \bg\ by defining on generators:
$$ C(s)(\langle t,u\rangle )=\langle st,u\rangle .
$$ The natural map $\delta:C\rightarrow \End_\bg$ has components
$\delta_x:C(x)\rightarrow \End_\bg(x)$ which are defined on generators as:
$$
\delta_x(\langle t,u\rangle )=tf(u)t^{-1}.
$$ It is easy to show that this defines $F$ on objects. Now on arrows: For an
arrow $(\alpha,\beta):(X,f,\bg)\rightarrow(X',f',\bg')$ in
\agpd, we define $F(\alpha,\beta)=(\beta,\overline{\alpha})$, where
$\overline{\alpha}:C\rightarrow C'\circ \beta$ has components defined on
generators by
$$
\overline{\alpha}_x(\langle t,u\rangle )= \langle \beta(t),\alpha(u)\rangle .
$$ In this way we get a functor $F:\agpd\rightarrow\pxm$. This is easily
verified to be left adjoint to $U'$ (see \cite{Garcia2003}{, Proposici\'on
3.1.13} for the details).
\end{proof}

\begin{proposition}
\label{pxm tripleable over agpd} The composite functor
$U_2:\xm\xrightarrow{U}\pxm\xrightarrow{U'}\agpd$ is monadic.
\end{proposition}

\begin{proof} We already know that $U_2$ has a left adjoint. By Beck's
tripleability theorem it is sufficient to prove that it reflects isomorphisms
and that it preserves coequalizers of $U_2$-contractible pairs. The first thing
is easy to see. The second thing requires a careful analysis of coequalizers in
$\xm$ and in $\agpd$. It is proved in this way: on the coequalizer of the
$U_2$-image of a $U_2$-contractible pair one can build in a natural way a
structure of crossed module together with a map in $\xm$ from the codomain of
the given contractible pair to this crossed module. After some tedious
calculations one verifies that this map is a coequalizer in $\xm$ and that its
image is the coequalizer in $\agpd$ from which it was built, proving the desired
property. See \cite{Garcia2003}, Proposici\'on 3.1.15, p. 148, for the details.
\end{proof}

We denote $\cotrip[2]$ the cotriple induced on $\xm$ by the the monadic
functor $U_2$, that is, $\cotrip[2] = F_2U_2$.

The next proposition reminds us of the well known fact that crossed modules can
be regarded as groupoids (actually, as 2-groupoids or groupoids enriched in the
category of groupoids). For some purposes in our context it will be convenient to
regard 2-groupoids (i.e. crossed modules) as those special double groupoids
(or internal groupoids in the category of groupoids) whose groupoids of objects and
of arrows have the same set of objects and whose structural functors (domain,
codomain, identity, and composition) are the identity on objects.

\begin{proposition}
\label{xm as groupoids} There is a functor  $\fxm:\Gpd(\Gpd)\rightarrow\xm$,
from the category of double groupoids to that of crossed modules, which has a
pseudo section
\begin{equation} \label{eq xm as gpd}
\gpd:\xm\rightarrow \Gpd(\Gpd),
\end{equation}
allowing us to regard any crossed module $(\bg,C,\delta)$ as an internal groupoid
in groupoids, having \bg\ as groupoid of objects and with groupoid of arrows given
by the ``semidirect product'' or Grothendieck construction $\semi{\bg}{C}=\int_\bg
C$. Furthermore, the above functors establish an isomorphism between the category
of crossed modules and the category of 2-groupoids.
\end{proposition}

\noindent(See Proposition \ref{isotilden} for a ``higher dimensional version'' of
\ref{xm as groupoids}.)

\begin{proof} Let $(\bg_0,\bg_1,s,t,i)$ be the underlying reflexive graph of a
double groupoid \cbg. We define a pre-crossed module $\fxm(\cbg)=(C,\delta)$ over
the groupoid $\bg_0$ of objects of \cbg\ by $$C=(\ker\widetilde{s})\circ i \quad
\text{and} \quad \delta=(\widetilde{t}\circ j)\ast i,$$ where
$j:\ker\widetilde{s}\rightarrow \End_{\bg_1}$ is the canonical inclusion, and we
use again the notation of tilde to denote the natural transformation
$\widetilde{f}:\End_\bg\rightarrow\End_{\bg'}\circ f$ induced by a
functor $f:\bg\rightarrow\bg'$. It is immediate to verify that $\delta$
satisfies Peiffer's identity and therefore $\fxm(\cbg)$ is a crossed module. This
defines the functor $\fxm$ on objects. On arrows $(f_0,f_1): \cbg \rightarrow
\cbg'$ $\fxm$ is defined by $\fxm(f_0,f_1) =
\textcolor{black}{(f_0,\alpha)}$, where $\alpha=(\widetilde{f_1}\circ j)\ast i$.

Let us now define the functor $\gpd:\xm\rightarrow\Gpd(\Gpd)$.  Given a
\bg-crossed module $\c=(C,\delta)$ by applying Grothendieck semidirect product
construction to $C$ we get a groupoid $\semi{\bg}{C}$ together with a canonical
split projection
$$
\xymatrix@C=1pc@R=1.75pc {**[l]\semi{\bg}{C}\ar@<-.2ex>[rr]_-s && \bg
\ar@/_2ex/[ll]_-i}.
$$
which is the identity on objects. Then the underlying reflexive graph of
$\gpd(\c)$ is (\bg, $\semi{\bg}{C}$, $s,t,i$), were the functor $t$ is the
identity on objects and takes any arrow $\textcolor{black}{(u,a)}:x\rightarrow
y$ in $\semi{\bg}{C}$ ($u:x\rightarrow y$ an arrow in $\bg$ and $a\in C(y)$) to
the composition $\delta_y(a)\circ u$. The composition map making this graph into
an internal groupoid in groupoids is the only possible one, which on the arrows
$x\to y$ of $\semi{\bg}{C}$ is given by the formula
$$
\textcolor{black}{ (v,b)\circ (u,a)=(u,b a) \qquad \bpar{\text{supposed} \ v =
\delta_y(a)u}. }
$$ This double groupoid is in fact a 2-groupoid. To an arrow $(f,\alpha) :
(\bg,C,\delta) \to (\bg',C',\delta')$  $\gpd$ associates the map of crossed
modules $\gpd(f,\alpha) = (f,\alpha')$ where $\alpha':\semi\bg C \to
\semi{\bg'}{C'}$ is the functor defined by $$\alpha'(u,a) =
\bpar{f(u),\alpha_y(a)}$$ for each $u:x\to y$ and $a\in C(y)$. It is easy to
verify that the crossed module corresponding to $\gpd(\c)$ is  isomorphic to \c\
and also that for any 2-groupoid \cbg\ the 2-groupoid $\gpd
\big(\fxm(\cbg)\big)$ is isomorphic to \cbg.
\end{proof}

Note that for any groupoid $\bg$ the functors $\fxm$ and $\gpd$ induce an
equivalence between the category $\xm_\bg$ and the subcategory of those
2-groupoids determined by those 2-groupoids having $\bg$ as groupoid of
objects and by those functors which are the identity on objects. \label{gpd of
conn comps of a xm} We note also that, the fundamental groupoid of a crossed
module \c\ is equal to the groupoid of connected components of the 2-groupoid
$\gpd(\c)$.

\textcolor{black}{By the set of connected components of a (pre-) crossed module,
$\pi_0(\c)$, it is understood the set of connected components of its base
groupoid.} Besides $\pi_0(\c)$ and $\pi_1(\c)$, the commutativity of diagram
\eqref{consequence xm} allows us to define the second ``homotopy group'',
$\pi_2(\c)$, of the crossed module $\c$, as the unique $\pi_1(\c)$-module such
that $\pi_2(\c)\circ q=\ker\delta$,
\begin{equation}
\label{first kernel induced map} \vcenter{ \xymatrix@C=1.5pc@R=1.75pc@!=1em {
{\vphantom{\Big(}}\widehat{C}
       \ar@<.6ex>[r]^-{\widehat{\delta}}
       \ar@<-.4ex>[r]_-0 & {\bg}
       \ar@/_.1 pt/[dr]_{\ker\delta}
       \ar[rr]^-q &  & **[r]{\pi_1(\c).}
       \ar@/^.1 pt/@{.>}[dl]^{\pi_2(\c)}
\\  &  &
\Ab & } }
\end{equation}

\section{Crossed complexes and their Postnikov towers}
\label{crs}

As indicated in the Introduction, we give a definition of crossed complex which
rests on the concept of crossed module instead of (the usual) building crossed
complexes all the way from groupoids. Crossed complexes over a fixed groupoid
are very easy to define as special types of chain complexes in the category of
crossed modules over the given groupoid. Having done that, it is evident how to
define morphisms between crossed complexes over different groupoids to get the
full category of crossed complexes. Standard references for crossed complexes
are \cite{BrHi1981a}, \cite{BrHi1981b}, and \cite{Tonks1993}.

If \bg\ is a fixed groupoid, a chain of complex in $\xm_\bg$ is a diagram
$$
\xymatrix@C=1.75pc@R=1.75pc { \cdots \ar[r] & \c_{n+1}\ar[r]^-{\partial_{n+1}} &
\c_n\ar[r]^-{\partial_n} &
\c_{n-1}\ar[r] & \cdots \ar[r] & \c_2\ar[r]^-{\partial_2} & \c_1,}
$$ of \bg-crossed modules whose underlying diagram of \bg-groups is a chain
complex in $\Gp^\bg$. In such a diagram, the fact that (for $n>1$) there is a
zero map from $\c_{n+1}$ to $\c_{n-1}$, in $\xm_\bg$, implies that in the
crossed module $\c_{n+1}=(C_{n+1},\delta_{n+1})$, $\delta_{n+1}=0$ and therefore
$\c_{n+1}$ is abelian, meaning that it is not just a \bg-group but a \bg-module,
$C_{n+1}:\bg\rightarrow\Ab$. Thus, in a chain complex in $\xm_\bg$ such as the
above one, for all $n\geq3$, $\c_n$ is an abelian crossed module.

\begin{definition} If \bg\ is a groupoid, a \bg-crossed complex is a chain
complex in $\xm_\bg$ of the form
\begin{equation}
\label{a crs}
\cbc:
\xymatrix@C=1.75pc@R=1.75pc { \cdots \ar[r] & \c_{n+1}
\ar[r]^-{\partial_{n+1}} & \c_n\ar[r]^-{\partial_n} &
\c_{n-1}\ar[r] & \cdots \ar[r] &
\c_2\ar[r]^-{\partial_2} & \termxm[\bg],}
\end{equation} such that for $n\geq 3$ the action of $\Img(\partial_2)$ on
$\widehat{C}_n$ is trivial. In other words, for every $n\geq3$ the following
diagram commutes
\begin{equation}
\label{condition crs}
\everyentry={\vphantom{\Big(}}
\xymatrix@C=1.75pc@R=1.75pc
{\widehat{C}_2\ar@<.6ex>[r]^{\widehat{\delta}_2}\ar@<-.4ex>[r]_0 & \bg
\ar[r]^-{C_n} &
\Ab.}
\end{equation} Here we are using the notation $\c_n=(C_n,\delta_n),\; n\geq 2,$
for the crossed modules in \cbc. The groupoid \bg\ is called the base groupoid
of the crossed complex, and $\c_2$ is called the base crossed module.
\end{definition}

For $n\geq3$, by the commutativity of \eqref{condition crs}, $C_n$ induces a
$\pi_1(\c_2)$-module $\overline{C}_n$,
\begin{equation}
\label{induced module by a ch complex} \vcenter{ \xymatrix@C=1.6pc@R=1.75pc@!=1em {
{\vphantom{\Big(}}\widehat{C}_2
       \ar@<.6ex>[r]^-{\widehat{\delta}_2}
       \ar@<-.4ex>[r]_-0 & {\bg}
       \ar@/_.1 pt/[dr]_{C_n}
       \ar[rr]^-q & & **[r]{\pi_1(\c_2)}
       \ar@/^.1 pt/@{.>}[dl]^{\overline{C}_n}
\\ & &
\Ab & } }
\end{equation} and the crossed complex \eqref{a crs} induces a chain complex of
$\pi_1(\c_2)$-modules of the form,
\begin{equation}\label{induced chain complex}
\overline{\cbc}:
\xymatrix@C=1.5pc@R=1.5pc {\cdots \ar[r] &
\overline{C}_{n+1}\ar[r]^-{\partial_{n+1}} &
\overline{C}_n\ar[r]^-{\partial_{n}} &
\overline{C}_{n-1}\ar[r] & \cdots \ar[r] & \overline{C}_3\ar[r]^-{\partial_{3}}
& \pi_2(\c_2)\ar[r]^-0 & 0,}
\end{equation} which will be used in the definition of the higher ``homotopy
groups'' of a crossed complex. (We name the natural maps in this chain complex
after their corresponding maps in \eqref{a crs} because they are essentially the
same, having the same components in $\Ab$.)

If \cbc\ is a \bg-crossed complex and $\cbc'$ is a $\bg'$-crossed
complex, a morphism $\bf: \cbc \rightarrow \cbc'$ is just a chain map,
that is, a family $\bf = \{f_n: \c_n \rightarrow \c'_n\}_{n\geq 1}$ of
maps of crossed modules such that for $n\geq 1$, $f_n\partial_{n+1}
=\partial'_{n+1}f_n$,
$$
\xymatrix{
 \c_{n+1} \ar[d]_{f_{n+1}} \ar[r]^{\partial_{n+1}} & \c_n \ar[d]^{f_n} \\
 \c'_{n+1} \ar[r]^{\partial'_{n+1}} & \c'_n .   }
$$ Note that this condition implies that all maps $f_n$ have the same
change-of-base functor, which is equal to $f_1$. The resulting category of
crossed complexes will be denoted \crs.

A morphism of crossed complexes $\bf:\cbc\rightarrow\cbc'$ is a
\emph{fibration} if each component $f_n$, is a fibration of crossed modules,
that is, if the functor $f_1:\bg\rightarrow\bg'$ is a fibration of
groupoids and the natural map of \bg-groups underlying  each $f_n$ is
surjective. The fibrations in \crs\ are part of a Quillen model structure in
this category.

 \gdef\thefundcrs{ The ``fundamental crossed complex" functor (see
\cite{BrHi1991}) associates a crossed complex to each simplicial set, $$
\Pi:\sset\rightarrow\crs. $$ This functor has a right adjoint ``nerve of a
crossed complex" which associates to each crossed complex \cbc\ the simplicial
set whose set of $n$-simplices is the set of maps of crossed complexes $$
\ner(\cbc)_n=\crs\bpar{\Pi(\Delta[n]),\cbc}. $$ Applying this nerve functor to a
fibration of crossed complexes one obtains a fibration of simplicial sets whose
fibers have the same homotopy type as those of the initial fibration. \vskip0pt
Since the geometric realization functor $\sset\rightarrow\Top$ also preserves
fibrations and the homotopy type of their fibers, we can define a functor
``classifying space of a crossed complex" which again has the same property. If
we restrict ourselves to the category \t\ of those spaces having the homotopy
type of a crossed complex we obtain an adjoint pair \begin{equation}\label{fund
crs ladj to geom realiz} \Pi\dashv B,\; \adj(\Pi,\t,\crs,B) \end{equation} which
induces an equivalence in the corresponding homotopy categories and which allow
us to reduce the calculation of the Postnikov towers of the spaces in \t, to the
calculation of Postnikov towers in \crs.
 }

An $n$-crossed complex or crossed complex of rank $n$ is a crossed complex such
as \eqref{a crs} in which all crossed modules $\c_m$ for $m>n$ are equal to
$\initxm[\bg]$. The full subcategory of \crs\ determined by the $n$-crossed
complexes will be denoted $\crs_n$. For a \bg-crossed complex to be of rank $0$
it is necessary that \bg\ be a discrete groupoid, that is, just a set.
Conversely, associated to a discrete groupoid \bg\ there is precisely one
$0$-crossed complex over \bg. Thus, $\crs_0$ can be identified with the category
of sets and we will put $\crs_0=\sets$. Similarly, since a map between two
1-crossed complexes is completely determined by the change-of-base functor,
which may be arbitrary, we will identify $\crs_1$ with the category of groupoids
and we write $\crs_1=\Gpd$. Finally, we will also write $\crs_2=\xm$ for similar
reasons.

The objects in $\crs_n$ are homotopy $n$-types, that is, they have trivial
``homotopy groups" in dimensions greater than $n$. The homotopy groups of a
crossed complex are defined as follows: $\pi_0(\cbc)$ is the set of connected
components of the base groupoid, so $\pi_0(\cbc)=\pi_0(\bg)$. Similarly,
$\pi_1(\cbc)=\pi_1(\c_2)=\bg/\Img(\delta_2)$, the fundamental groupoid of the
base crossed module of \cbc. For $n\geq 2$, $\pi_n(\cbc)$ is defined as the
``homology group" $H_n(\overline{\cbc}):\pi_1(\cbc)\rightarrow \Ab$ of the
induced chain complex of $\pi_1(\cbc)$-modules \eqref{induced chain complex}.
Note that if we consider $\pi_n(\cbc)$ as a \bg-module via the canonical
projection $q:\bg\rightarrow\pi_1(\cbc)$, for $n\geq 2$, the \bg-crossed module
$(\pi_n(\cbc),0)$ is the kernel of the induced map
$\overline{\partial}_n:\c_n/\Img(\partial_{n+1})\rightarrow \c_{n-1}$ (see
\eqref{ind part} below).

\bigskip

In the same way that the discrete inclusion of sets into groupoids is both
reflexive and coreflexive, Proposition \ref{adjs betw xm and gpd} tells us that
the ``discrete" inclusion $\bg\mapsto \termxm[\bg]$ of groupoids into crossed
modules is both reflexive and coreflexive. These are particular cases of a
general situation. For every $n\geq 0$, the subcategory $\crs_n$ of $\crs$ is
both reflexive and coreflexive. We are mainly interested in the reflector
$\tilde{P}_n:\crs\rightarrow\crs_n$, left adjoint to the inclusion
$i_n:\crs_n\rightarrow\crs$. For $n=0,1,$ we have $\tilde{P}_0=\pi_0\circ \base$
(``set of connected components of the base groupoid") and
$\tilde{P}_1=\pi_1\circ\base$ (``fundamental groupoid of the base crossed
module"). For higher $n, \; \tilde{P}_n$ is calculated in terms of the following
coequalizer in \xm:
$$
\xymatrix@C=1.75pc@R=1.75pc
{\c_{n+1}\ar@<.6ex>[r]^-{\partial_{n+1}}\ar@<-.4ex>[r]_-0 &
\c_n \ar[r]^-{q_n} &
\c_n/\Img\partial_{n+1}. }
$$ Thus $\tilde{P}_n$ associates to the crossed complex \cbc\ given in \eqref{a
crs} the following $n$-crossed complex:
\begin{equation}
\label{ind part}
\tilde{P}_n(\cbc):
\xymatrix@C=1.25pc@R=1.75pc { \cdots \ar[r] & \initxm[\bg] \ar[r] &
\c_n/\Img\partial_{n+1}
\ar[r]^-{\overline{\partial}_n} & \c_{n-1}\ar[r] & \cdots \ar[r] &
\c_2\ar[r]^-{\partial_2} &
\termxm[\bg] \; ,}
\end{equation} where $\overline{\partial}_n$ is the unique map of crossed
modules such that $\partial_n=\overline{\partial}_n\circ q_n$, induced by the
fact that $\partial_n\partial_{n+1}=0$.

The objects in $\crs_n$ are homotopy $n$-types and, in addition, all $n$-types
of crossed complexes are represented in $\crs_n$. That is, if $\cbc\in \crs$ is
any crossed complex which is an $n$-type, there exists an object
$\cbc_n\in\crs_n$ which is homotopically equivalent to \cbc. That object is just
$\cbc_n=\tilde{P}_n(\cbc)$.

Regarding the reflectors $\tilde{P}_n$ as endofunctors, $P_n$, of \crs, we have
a situation as described in the introduction. We have idempotent endofunctors
$P_n:\crs\rightarrow\crs$ such that $P_n=P_{n+1}P_n$, and
$\eta_{n+1}=P_{n+1}*\delta^{(n)}$ is the composition of $P_{n+1}$ with the unit
$\delta^{(n)}$ of the adjunction $\tilde{P}_n\dashv i_n$. Working out the
components of the unit $\delta^{(n)}$, for $n>1$, one finds that for a given
\bg-crossed complex \cbc, the components of the map
$(\eta_{n+1})_\cbc:P_{n+1}(\cbc)\rightarrow P_n(\cbc)$ are: the trivial map
$\c_{n+1}/\Img\partial_{n+2}\rightarrow \initxm[\bg]$ at dimension $n+1$, the
``projection to the quotient", $q_n:\c_n\rightarrow \c_n/\Img\partial_{n+1}$, at
dimension $n$, and an identity map at all other dimensions. The cases of
$\eta_1$ and $\eta_2$ are little different since for these maps the
change-of-base functor is not an identity. For $\eta_1$ the change of base
functor is the canonical projection $q_0:\pi_1(\cbc)\rightarrow\pi_0(\cbc)$,
while for $\eta_2$ it is the canonical projection
$q_1:\bg\rightarrow\pi_1(\cbc)$.

\begin{proposition}
\label{fibras son kpin} For every crossed complex $\cbc\in\crs$, and every
$n\geq0$ the map $(\eta_{n+1})_\cbc:P_{n+1}(\cbc)\rightarrow P_n(\cbc)$ is a
fibration with fibers of the type of $K(\Pi,n+1)$. For $n>1$, the fiber of
$(\eta_{n+1})_\cbc$ over $x\in\bg$ has the homotopy type of
$K\bpar{\pi_{n+1}(\cbc)(x),n+1}$.
\end{proposition}
\begin{proof} Let us first consider $(\eta_1)_{\cbc}$, which is
$q_0:\pi_1(\cbc)\rightarrow \pi_0(\cbc)$, a surjective map to a discrete
groupoid and therefore it is a fibration of groupoids. The fiber over a given
connected component $\bar{x}\in\pi_0(\cbc)$ is a connected groupoid and
therefore has the homotopy type of a $K(\Pi,1)$ (taking for $\Pi$ any of the
groups of endomorphisms of any object in that connected groupoid).

Next, we look at $(\eta_2)_\cbc:P_2(\cbc)\rightarrow P_1(\cbc)$. The fiber of
this map over an object $x\in\bg$ is the reduced 2-crossed complex
$$ (\c_2/\Img\partial_3)(x):C_2(x)/\Img(\partial_3)_x\rightarrow
\Img(\partial_2)_x\; .
$$ This is a crossed module over the group $\Img(\partial_2)_x$ and therefore it
has $\pi_0=0$. Since the above map is surjective, this crossed module has
$\pi_1=0$, while $\pi_2$ is precisely the abelian group $\pi_2(\cbc)(x)$. For
higher $n$, $\pi_n=0$, thus $\eta_2$ is a fibration with fiber over $x$ of the
type $K\bpar{\pi_1(\cbc)(x),2}$.

For $n>2$ all the $\eta_n$ are morphisms of crossed complexes whose change of
base functor is the identity on objects (as in the case $n=2$). Therefore, their
fiber on an object $x\in\bg$ is a reduced crossed complex. The special thing for
$n>2$ is that the base groupoid of the fiber is trivial and therefore the fiber
is just a chain complex of abelian groups. In general, the fiber of
$(\eta_n)_\cbc$ (for $n>2$) over an $x\in\bg$ is a crossed complex with trivial
components below the $n-1$ an a surjective morphism in dimension $n$. Therefore,
all homotopy groups of the fiber are trivial in dimensions other than $n$, and
it is equal to $\pi_n(\cbc)(x)$ in dimension $n$.
\end{proof}

\begin{proposition}
\label{prop ptow of a crs} For every crossed complex \cbc, the chain of
fibrations
\begin{equation}
\label{diagr ptow of a crs}
\xymatrix@C=1.75pc@R=1.75pc { \cdots \ar[r]^-{\eta_{n+2}} &
P_{n+1}(\cbc)\ar[r]^-{\eta_{n+1}}& P_n(\cbc)\ar[r]^-{\eta_n} &\cdots
\ar[r]^-{\eta_1} & P_0(\cbc)}
\end{equation} is a Postnikov tower for \cbc.
\end{proposition}
\begin{proof} By Proposition \ref{fibras son kpin} it is sufficient to prove
that the limit of diagram \eqref{diagr ptow of a crs} is \cbc. We know that the
morphisms $\delta^{(n)}:\cbc\rightarrow P_n(\cbc)$ determined by the units of
the adjunctions $\tilde{P}_n\dashv i_n$ constitute a cone over \eqref{diagr ptow
of a crs}. Given any other cone $\{\phi^{(n)}:\cbc'\rightarrow
P_n(\cbc)\}$ over \eqref{diagr ptow of a crs} there is a unique way of defining
a map of crossed complexes $\bf:\cbc'\rightarrow\cbc$ such that for all
$n$, $\delta^{(n)}\bf=\phi^{(n)}$. One just needs to take into account that
$\delta^{(n)}_m$ is an identity map for all $m<n$ and define
$f_n=\phi_n^{(n+1)}$.
\end{proof}

The subcategories $\crs_n$ of \crs\ are not only reflexive, but also
coreflexive, the right adjoint to the inclusion being ``simple truncation",
$T_n:\crs\rightarrow\crs_n$, so that $T_0$ is essentially the set of objects of
the base groupoid, $T_1$ is ``base groupoid", $T_2$ is ``base crossed module".
Furthermore, $T_n$ has itself a right adjoint denoted $\cosk^n$. For $n=0$ and
$n=1$ this further right adjoint is ``codiscrete groupoid on a set" ($n=0$) and
``trivial crossed module on a groupoid" \bpar{so that
$\cosk^1(\bg)=\termxm[\bg]$ or $\cosk^1(\bg)=(\cdots
\initxm[\bg]\rightarrow\cdots\rightarrow \initxm[\bg]\rightarrow
\termxm[\bg])$}. For $n>1$, the right adjoint to $T_n$,
$\cosk^n:\crs_n\rightarrow\crs$, assigns to an $n$-crossed complex
$$
\cbc:
\xymatrix@C=1.25pc@R=1.75pc { \cdots\ar[r] & \initxm[\bg] \ar[r] &
\initxm[\bg]\ar[r] & \c_n\ar[r]^-{\partial_n} &
\c_{n-1}\ar[r] & \cdots \ar[r] & \c_2\ar[r]^-{\partial_2} & \termxm[\bg],}
$$ the following $(n+1)$-crossed complex:
$$
\cosk^n(\cbc):
\xymatrix@C=1.25pc@R=1.75pc { \cdots\ar[r] & \initxm[\bg] \ar[r] &
\ker\partial_n\ar[r] &
\c_n\ar[r]^-{\partial_n} & \c_{n-1}\ar[r] & \cdots \ar[r] &
\c_2\ar[r]^-{\partial_2} & \termxm[\bg]\; .}
$$ For $n>2$, the functor \begin{equation}
\label{ntecho}\ntecho:\crs_{n,\bg}\rightarrow\ab^\bg \end{equation} takes each $n$-crossed
complex \cbc\ having \bg\ as  base groupoid, to the \bg-module
$\ntecho(\cbc)=\techo(\c_n)$. Note that, since finite limits and
coequalizers in $\crs_{n,\bg}$ are calculated componentwise, and since
the functor $\techo = \techo_2$ preserves finite limits and
coequalizers, for $n>2$ the functor \ntecho\ also preserves finite
limits and coequalizers.

We end this section with a higher dimensional analog of Proposition \ref{xm as
groupoids}. The idea is to regard the $(n+1)$-crossed complexes as some kind of
internal groupoids  in the category of $n$-crossed complexes. For $n>1$ there is
a difficulty we do not have in Proposition \ref{xm as groupoids}. For example,
(in case $n=2$) it is possible to carry out a construction analogous to the one
defining the functor $\fxm$, but starting with a groupoid internal in crossed
modules: $(\c_0,\c_1,s,t,i,\gamma)$. Let the crossed modules of objects and
arrows of this groupoid be $\c_i=(\bg_i,C_i,\delta_i)$, $i=0,1$, and let the
domain, codomain and identity maps be $s=(f_s,\alpha_s)$, etc. We can define
$\c_3\xto{\partial_3}\c_2\xto{\delta}\termxm[\bg_0]$ where
$\c_3=(\bg_0,\ker(\alpha_s*f_i),\delta\partial_3)$, $\delta=\delta_0$,
$\partial_3=\alpha_t*f_i$, and $\c_2=\c_0$. One does not, however, obtain
directly from this construction a 3-crossed complex unless the base groupoids of
the crossed modules $\c_0$, $\c_1$ are the same and the structural maps ($s$,
$t$, $i$ and $\gamma$) have trivial change of base (i.e. $f_s=1_{\bg_0}$ etc.).
It is of course possible to force the result to be a 3-crossed complex by making
the appropriate quotients, but this would introduce unnecessary complication.
Thus, we shall appropriately restrict the categories of internal groupoids in
crossed complexes so that, for example, in the case $n=2$ we will only consider
those internal groupoids in crossed modules satisfying the conditions said
above. In general, for each $n>0$ let $\gcrs[n]$ be the full subcategory of
$\Gpd(\crs_n)$ (internal groupoids in $\crs_n$) determined by those groupoids
$\cbg\in\Gpd(\crs_n)$ whose $n$-crossed complex of objects has the same
$(n-1)$-truncation as its $n$-crossed complex of arrows and whose structural
maps (domain, codomain, identity, and composition) have the identity map as
$(n-1)$-truncation. Thus, an object $\cbg\in\gcrs[n]$ gives rise to a diagram in
$\xm$ of the form
\begin{equation}
\label{gcrsn}
\xy
   (0,0)*+{\c_n^1 \times_{\c_n^0} \c_n^1}="c",
   (30,0)*+{\c_n^1}="d",
   (60,0)*+{\c_n^0}="e",
   (30,-15)*+{\c_{n-1}}="g",
   (30,-25)*+{\c_{n-2}}="h",
   (30,-35)*+{\vdots}="i",
   (30,-45)*+{\c_2}="j",
   (30,-55)*+{\mathbf{1}_\g}="k",
\ar @/^-4ex/ "e";"d" |{id}
\ar @{->}^{s} "d";"e" <3pt>
\ar @{->}_{t} "d";"e" <-3pt>
\ar @{->}^-{\circ} "c";"d"
\ar @/_3ex/ "c";"g"_{\partial_n^1 \times_{\partial_n^0}\partial_n^1}
\ar @{->}^{\partial_n^1} "d";"g"
\ar @/^3ex/ "e";"g"^{\partial_n^0}
\ar @{->}^{\partial_{n-1}} "g";"h"
\ar @{->} "h";"i"
\ar @{->} "i";"j"
\ar @{->}^{\partial_2} "j";"k"
\endxy
\end{equation}

We can now state:

\begin{proposition}
\label{isotilden}  For each $n>1$ there is a functor  $\fcrs_n:\gcrs[n]\to
\crs_{n+1}$, which has a  section
\begin{equation}
\label{eq xm as gpd general n}
\gpd_n:\crs_{n+1}\to \gcrs[n],
\end{equation} allowing us to regard any $(n+1)$-crossed complex
$$\cbc:
\xymatrix@C=1.25pc@R=1.75pc { \cdots\ar[r] & \initxm[\bg]\ar[r] &
\c_{n+1}\ar[r]^-{\partial_{n+1}} & \c_n\ar[r]^-{\partial_n} &
\c_{n-1}\ar[r] & \cdots \ar[r] & \c_2\ar[r]^-{\partial_2} & \termxm[\bg],}
$$  as an internal groupoid in $n$-crossed  complexes, having
$\cbf_0=T_{n}(\cbc)$ as $n$-crossed complex of objects and whose $n$-crossed
complex of arrows is given by
$$
\cbf_1: \c_{n+1} \times \c_n \xto{\partial_n \, p_0} \c_{n-1}
\xto{\partial_{n-1}} \ldots \to \c_2 \to \mathbf{1}_\bg.
$$ where $\c_{n+1} \times \c_n$ is a cartesian product in $\xm_\bg$ and
$p_0:\c_{n+1}\times\c_n\to \c_n$ is the corresponding canonical projection. The
functors $\fcrs_n$ and $\gpd_n$ establish an isomorphism between $\crs_{n+1}$
and $\gcrs[n]$.
\end{proposition}

\begin{proof} Let us complete the definition of $\gpd_n$. The domain map
$s:\cbf_1 \to \cbf_0$ is induced by the projection $p_0:\c_{n+1} \times \c_n \to
\c_n$, and the codomain map $t:\cbf_1 \to \cbf_0$ is induced by the morphism of
\bg-groups $C_{n+1} \times C_n \to C_n$ (actually $\bg$-modules except in the
case $n=2$) defined on $(u,v)\in C_{n+1}(x) \times C_n(x)$ as
$(u,v)\mapsto\partial_{n+1}(u) \,v\in C_n(x)$ (note that even in the case $n=2$,
in which $C_2(x)$ may not be abelian, $\partial_{n+1}$ takes its values in the
center of $C_n$ and thus this is a homomorphism). Clearly the canonical map $C_n
\hookrightarrow C_{n+1} \times C_n$ is a common section for $s$ and $t$, so that
we get an internal graph
\begin{equation}
\label{grcrnt}\xymatrix{\cbf_1 \ar@<0.6ex>[r]^s \ar@<-0.6ex>[r]_t & \cbf_0
\ar@/_1.2pc/[l]_{id}}
\end{equation} in the category $(\crs_n)_{T_{n-1}(\cbc)}$ of $n$-crossed
complexes whose $(n-1)$-truncation is $T_{n-1}(\cbc)$.

We want to endow this graph with a structure of internal groupoid  in
$\crs_{n}$. For this it is sufficient to do it at the highest dimension, only
place where the graph structure is not trivial. In this dimension we have the
internal graph
\begin{equation}
\label{grab} \xymatrix{C_{n+1} \times C_n \ar@<0.6ex>[r]^-s \ar@<-0.6ex>[r]_-t &
C_n
\ar@/_1.5pc/[l]_{id}}
\end{equation} in the category $\Gp^\bg$ of \bg-groups, which admits a unique
structure of internal groupoid in $\Gp^\bg$.

The structure of internal groupoid in $\gp^\bg$ of the graph \eqref{grab}
determines a structure of internal groupoid in $(\crs_n)_{T_{n-1}(\cbc)}$ on the
graph \eqref{grcrnt},  this internal groupoid in  $\crs_n$ will be denoted
$\gpd_n(\cbc)$. It is easy to show that this construction is functorial and thus
we have a functor
\begin{equation}
\label{ngd} \gpd_n: \crs_{n+1} \to \gcrs[n].
\end{equation}
In order to define its cuasi-inverse $\fcrs_n:\gcrs[n]
\longrightarrow\crs_{n+1}$ we consider an object
$$
\cbg: \xymatrix{\cbc^1 \ar@<0.6ex>[r]^s \ar@<-0.6ex>[r]_t & \cbc^0
\ar@/_1.2pc/[l]_{id}}
$$ in \gcrs[n], as in \eqref{gcrsn}, and we apply \ntecho\ to the morphism
$s:\cbc^1\to \cbc^0$. We then obtain a morphism of \bg-groups
$$ s=\ntecho(s):\ntecho(\cbc^1)=C_n^1 \longrightarrow C_n^0=\ntecho(\cbc^0).
$$ whose kernel is  a $\bg$-module $K=\Ker(s):\bg \to \ab$ associating with each
object $x\in \bg$ the subgroup $K(x)$ of $C_n^1(x)$ consisting of the elements
$u \in C_n^1(x)$ such that $s_x(u)= 0_{C_n^0(x)}$, with an action which is
induced by the action of $C_n^1$. Evidently, this $\bg$-module determines a
crossed module $\ck=(\bg,K,0)$ and  $\partial_2$ acts trivially on the totally
disconnected groupoid $\widehat{K}$. As a result we have a $(n+1)$-crossed
complex
$$
\fcrs_n(\cbg)=(\ck \xto{\partial_{n+1}} \c_n^0 \xto{\partial^0_n} \c_{n-1}
\to \ldots \to \c_2 \xto{\partial_2} \mathbf{1}_\bg)
$$ where $\partial_{n+1}:K \to C^0_n$ is a morphism of $\bg$-groups induced by
the morphism of $\bg$-groups $t=\ntecho(t):C_n^1 \to C_n^0$ associated to the codomain of $\cbg$,
that is, for each object $x \in \bg$, the $x$-component of $\partial_{n+1}$ is given by
$$ (\partial_{n+1})_x: K(x) \to C^0_n(x), \; (\partial_{n+1})_x(u) = t_x(u).
$$ Note that  $\fcrs_n(\cbg)$ is really a chain complex, that is, $ \partial^0_n
\partial_{n+1} = 0.$ This construction of $\fcrs_n(\cbg)$ is also functorial so
that we have a functor \begin{equation} \label{fcrsn} \fcrs_n: \gcrs[n]  \to
\crs_{n+1}. \end{equation} Let us see that it is a cuasi-inverse for $\gpd_n$.

If $\cbg \in \gcrs[n]$ as in \eqref{gcrsn}, $\ntecho(\cbg)$ is an internal groupoid in the
category of $\bg$-groups, hence for every $x \in \obj(\bg)$ we have an internal groupoid in the
category of groups
$$
\xymatrix{C^1_n(x) \times_{C_n^0(x)} C^1_n(x) \ar[r] & C_n^1(x)
\ar@<0.6ex>[r]^{s_x}
\ar@<-0.6ex>[r]_{t_x} &  C_n^0(x) \ar@/_1.2pc/[l]_{id_x}}.
$$ Thus, we have a group isomorphism
$$ G_x:\semi{C_n^0(x)}{\N_{(s_x,id_x)}}=K(x) \times C_n^0(x) \xto{\cong}
C_n^1(x) ;\qquad G_x(u,v)= u \ id_x(v).
$$ Since the above isomorphism is natural, it is immediate that the pair $(G,
Id_{C_n^0})$ is an isomorphism of graphs in $\gp^\bg$ which induces a graph isomorphism in
$\crs_n$, hence an isomorphism in \gcrs[n] between the groupoids $\gpd_n\fcrs_n(\cbg)$ and
$\cbg$. Conversely, if \cbc\ is a $(n+1)$-crossed complex, then
$$
\fcrs_n\big(\gpd_n(\cbc)\big) = (\ck \xto{\partial_{n+1}} \c_n
\xto{\partial^0_n} \c_{n-1}
\to \ldots \to \c_2 \xto{\partial_2} \mathbf{1}_\bg),
$$ where $\ck = (\bg,K,0)$ with $K = \Ker(C_{n+1}\times C_n \xto{s}C_n)$. Since
$s$ is the canonical projection, it is clear that $K=C_{n+1}$ and $\ck = \c_{n+1}$.
Looking closely at the connecting morphisms in $\fcrs_n(\gpd_n(\cbc))$ one
realizes immediately that $\fcrs_n\big(\gpd_n(\cbc)\big) = \cbc$.
\end{proof}

\begin{proposition}
\label{nypn} For any $n>0$ the following equations of functors hold
$$
\xymatrix{\crs_{n+1}\ar[r]^{\incl_{n+1}}\ar[d]_{\gpd_n} & \crs
\ar[d]^{\widetilde{P}_n} \\
\gcrs[n] \ar[r]_{\pi_0} & \crs_n ,}\qquad
\xymatrix{\crs_{n+1}\ar[r]^{\incl_{n+1}} &\crs\ar[d]^{\widetilde{P}_n} \\
\gcrs[n]\ar[u]^{\fcrs_n}\ar[r]_{\pi_0} &\crs_n,}
$$
$$
\widetilde{P}_n\incl_{n+1} = \pi_0\gpd_n\,, \qquad
\widetilde{P}_n\incl_{n+1}\fcrs_n = \pi_0,
$$ in other words, for each $(n+1)$-crossed complex \cbc\ and each groupoid
$\cbg\in\gcrs[n]$:
$$ P_n(\cbc)=\pi_0\gpd_n(\cbc)\qquad  \mbox{and}   \qquad
\pi_0(\cbg)=P_n\fcrs_n(\cbg),
$$
\end{proposition}
\begin{proof} Since $\fcrs_n(\gpd_n(\cbc)) = \cbc$, it is sufficient to verify
the right hand equation, that is, $\pi_0(\cbg)=P_n\fcrs_n(\cbg)$. These two
crossed complexes agree in dimensions less than $n$. In dimension $n$,
$\pi_0(\cbg)$ is the coequalizer of $s$ and $t$ (see diagram \ref{gcrsn}). On
the other hand, $P_n\fcrs_n(\cbg)$ is, in dimension $n$, the quotient of $C^0_n$
by the image of $\partial_{n+1}(=t):\ker (s)\to C^0_n$. That this quotient is
equal to the previous coequalizer is an immediate consequence of the general
fact that the coequalizer of a parallel pair of group homomorphisms, $s,t:G\to
H$ having a common section is the quotient of $H$ by $t(\ker(s))$.
\end{proof}

\begin{proposition}
\label{endypi} For each groupoid $\cbg\in\gcrs[n]$ and each $(n+1)$-crossed complex \cbc\ we have
natural isomorphisms:
$$
\begin{array}{c}
\ntecho\End(\cbg)\cong \ntecho\obj(\cbg)\times (\pi_{n+1}\fcrs_n(\cbg)\circ \, q)
\qquad\mbox{y}\\ [2pc]
\ntecho\End\gpd_n(\cbc)\cong \ntecho T_n(\cbc)\times (\pi_{n+1}(\cbc)\circ\, q),
\end{array}
$$ where \bg\ denotes both the base groupoid of $\obj(\cbg)$ and that of \cbc,
and $q$ denotes either the canonical projection $\bg\to\pi_1\obj(\cbg)$ or
$\bg\to\pi_1(\cbc)$. \end{proposition} \begin{proof} The second isomorphism is a
consequence of the first one and of the identities $\obj(\gpd_n(\cbc))=T_n(\cbc)$
and $\fcrs_n \gpd_n(\cbc)=\cbc$.

For each object $x\in \bg$, the isomorphism
$$
\ntecho\End(\cbg)(x) \xto{\cong} \ntecho\obj(\cbg)(x)\times
(\pi_{n+1}\fcrs_n(\cbg)\circ
\, q(x))
$$ takes each $u\in \ntecho\End(\cbg)(x)$ to the pair
$$ (\; s(u)=t(u)\;,\; u-id(s(u))\; ) \in\ntecho\obj(\cbg)(x)\times
(\pi_{n+1}\fcrs_n(\cbg)\circ \, q(x)).
$$
\end{proof}

\section{The Postnikov invariants of a crossed complex} \label{inv}

As indicated in the Introduction, we distinguish two types of Postnikov
invariants of a crossed complex. On the one hand we have the ``\emph{algebraic}"
invariants, which are elements of \emph{algebraic} (cotriple) cohomologies in
the categories $\crs_n$. On the other hand, corresponding to each algebraic
invariant $k_{n+1}$, we have a ``\emph{topological}" invariant, which is an
element of a \emph{singular} cohomology. In this section we first determine the
algebraic invariants, characterizing them as extensions and as torsors. Then, we
define the topological invariants and the singular cohomologies in which they
live. Finally, we show how to map the cotriple cohomologies to the singular ones
so that one can obtain the topological invariants from the algebraic ones.

\subsection{The algebraic invariants}

Let \cbc\ be a $\bg$-crossed complex. For $n\geq 0$, the $(n+1)^{th}$ algebraic
Postnikov invariant, $k_{n+1}$, of \cbc\ is determined by the fibration
$\eta_{n+1}:P_{n+1}(\cbc)\rightarrow P_n(\cbc)$, which is completely described
by the following diagram:
\begin{equation} \label{eta n+1}
\vcenter{
\everyentry={\vphantom{\big(}}
\xymatrix@C=1.25pc@R=1.5pc {
\cdots\ar[r] & \initxm[\bg]\ar[r] \ar@{=}[d] & \c_{n+1} / \Img\partial_{n+2}
\ar[d]\ar[r]^-{\bar{\partial}_{n+1}} & \c_n\ar[r]^{\partial_n} \ar[d]^{q_n} &
\c_{n-1}\ar[r]^{\partial_{n-1}}\ar@{=}[d] & \cdots \\ \cdots \ar[r] &
\initxm[\bg]\ar[r] & \initxm[\bg]\ar[r] & \c_n/\Img\partial_{n+1}
\ar[r]^-{\bar{\partial}_n} & \c_{n-1}\ar[r]^-{\partial_{n-1}} & \cdots } }
\end{equation} For $n=0$ we get the first invariant, $k_1$, determined by
$\eta_1$:
$$
\everyentry={\vphantom{\big(}}
\xymatrix@C=1.25pc@R=1.5pc {\cdots \ar[r] &
\initxm[\pi_1(\cbc)]\ar@<-1.2ex>[d]^{(q_0,0)} \ar[r] & \termxm[\pi_1(\cbc)]
\ar@<-1.2ex>[d]^{(q_0,0)}
\\
\cdots \ar[r] & \initxm[\pi_0(\cbc)] \ar[r] & \termxm[{\pi_0(\cbc)}]}
$$ where we have indicated in the vertical arrows the change of base functor
$q_0$ since it is not an identity.

This invariant is an element in the topos cohomology of $\crs_0=\sets$, a
category which has very little structure and, correspondingly, with a very
simple cohomology: it is trivial in dimensions $\geq 1$, so that
$H^2(P_0(\cbc),\pi_1(\cbc))$ only has one element. Although it is not difficult
to see that this element corresponds to $\eta_1$, there is really no need of our
machinery to determine it, and we will not discuss here this case any further.
For a more complete discussion of this case we refer the interested reader to
\cite{Garcia2003}.

For $n=1$ we have the second invariant, $k_2$, determined by $\eta_2$:
$$
\everyentry={\vphantom{\big(}} \xymatrix@C=1.25pc@R=1.5pc {\cdots \ar[r] &
{\initxm[\bg]\ }\ar@<-1ex>[d]^{(q_1,0)} \ar[r] & \c_2/\Img\partial_3
\ar@<-1ex>[d]^{(q_1,0)}
\ar[r]^-{\bar{\partial}_2} & {\termxm[\bg]\ }\ar@<-1ex>[d]^{(q_1,q_1).} \\
\cdots \ar[r] &
\initxm[\pi_1(\cbc)] \ar[r] & \initxm[\pi_1(\cbc)]\ar[r] &
\termxm[{\pi_1(\cbc)}]}
$$ Note, however, that $q_1$ is the cokernel of $\bar{\partial}_2$, while
$\ker\bar{\partial}_2$ is, as noted earlier (see page \pageref{induced module by
a ch complex}), the \bg-crossed module $\a_2=(\pi_2(\cbc),0)$, where
$\pi_2(\cbc)$ is considered as \bg-module via $q_1:\bg\rightarrow\pi_1(\cbc)$.
Thus $\eta_2$ is completely determined by the following sequence of crossed
modules
\begin{equation}
\label{k2 ex seq xm} \everyentry={\vphantom{\big(}} \xymatrix@C=1.25pc@R=1.5pc {
0 \ar[r] & \a_2\ar[r] & \c_2/\Img\partial_3 \ar[r]^-{\bar{\partial}_2} &
\termxm[\bg] \ar[r]^-{q_1} & \termxm[{\pi_1(\cbc)}]\ar[r] & 0,}
\end{equation} which can be regarded as a genuine exact sequence
\begin{equation}
\label{k2 ex seq gpd} \everyentry={\vphantom{\Big(}} \xymatrix@C=1.25pc@R=1.5pc
{ 0 \ar[r] & \widehat{\pi_2(\cbc)} \ar[r] &
\widehat{C}_2/\Img\widehat{\partial}_3 \ar[r]^-{\bar{\partial}_2} & \bg
\ar[r]^-{q_1} & \pi_1(\cbc)\ar[r] & 0,}
\end{equation} in the category $\Gpd_{\obj(\bg)}$.

A sequence such as \eqref{k2 ex seq xm} or \eqref{k2 ex seq gpd} is an extension
of the groupoid $\pi_1(\cbc)$ by the $\pi_1(\cbc)$-module $\pi_2(\cbc)$,
according to the following definition, of which the ``\emph{reduced}" or pointed
case is the well known definition of 2-extensions of groups by modules:

\begin{definition}
\label{2exgpd} A 2-extension of a groupoid $\Pi$ by a $\Pi$-module
$A:\Pi\rightarrow\Ab$ is a crossed module $\c=(\bg,C,\delta)$, called the fiber
of the extension, together with an exact sequence

\begin{equation}
\everyentry={\vphantom{\Big(}}
\xymatrix@C=1.25pc@R=1.5pc { 0\ar[r] & \widehat{A} \ar[r] &
\widehat{C}\ar[r]^-{\widehat{\delta}} & \bg\ar[r]^-q & \Pi\ar[r] & 0}
\end{equation}

in $\Gpd_{\obj(\Pi)}$ such that the kernel of $\delta$ factors through
$\Pi$ as $A\circ q$,

\begin{equation}
\label{com tria}
\vcenter{
\xymatrix{\bg\ar[rr]^q\ar@/_.1pt/[dr]_-{\ker \delta} && \Pi \ar@/^.1pt/[dl]^A \\
& \Ab } }.
\end{equation}

\end{definition}

The standard definition of morphism of extensions gives rise to a category,
denoted
$\Ext^2(\Pi,A)$, whose objects are the 2-extensions of a groupoid $\Pi$ by  a
$\Pi$-module $A$. As usual we denote with brackets, $\Ext^2[\Pi,A]$, the
category of connected components of $\Ext^2(\Pi,A)$.

It is now easy to show that the extensions just defined represent cohomology
elements in a well known cotriple cohomology of groupoids. We use the notation
and results of the Appendix, Section~5.

\begin{proposition}
\label{2-exts of gpd are tor} Let $U:\Gpd\rightarrow\Gph$ be the underlying
graph functor defined on the category of groupoids. If $\Pi$ is a groupoid, $A$
is a $\Pi$-module, and $\abgc A1$ is the abelian group object in $\Gpd/\Pi$,
$\abgc A1 = (\semi\Pi A\leftrightarrows\Pi)$, there is a full and faithful
functor
$$
\Ext^2(\Pi,A)\rightarrow \Tor^2_U(\Pi,\tilde{A}_1).
$$
\end{proposition}

\begin{proof}
Consider a 2-extension of $\Pi$ by $A$, as in Definition \ref{2exgpd},
with
$\c=(\bg,C,\delta)$ being the fiber crossed module. The commutativity of the
triangle \eqref{com tria} gives a commutative square of groupoids and functors
$$
\everyentry={\vphantom{\Big(}}
\xymatrix@C=1.25pc@R=1.5pc {\semi{\bg}{\ker\delta}\ar[d]\ar[r]^-\alpha &
\semi{\Pi}{A} \ar[d] \\ \bg\ar[r]_-q &
\Pi,}
$$ which is a pullback. Let then $\tilde{\c}$ be the internal groupoid in
groupoids corresponding to the crossed module \c\ by the functor $\gpd$ of
Proposition \ref{xm as groupoids}. The groupoid of connected components of
$\tilde{\c}$ is $\pi_1(\c)=\Pi$, and the groupoid of endomorphisms of
$\tilde{\c}$ is $\semi{\bg}{\ker\delta}$. Therefore $(\Pi,\tilde{\c},\alpha)$ is
an $(\tilde{A}_1,2)$-torsor which is clearly $U$-split. Furthermore, it is a
routine straightforward verification to see that the above construction is
functorial. That the functor $\Ext^2(\Pi,A)\rightarrow
\Tor^2_U(\Pi,\tilde{A}_1)$ so defined is full and faithful is an immediate
consequence of Proposition \ref{xm as groupoids}.
\end{proof}

\begin{proposition}
\label{2-tor of gpd are ext} For any groupoid $\Pi$ and $\Pi$-module $A$,
$$ H^2_{\bbg_1}(\Pi,\tilde{A}_1)\cong \Ext^2[\Pi,A].
$$ where $\bbg_1$ is the cotriple on $\Gpd/\Pi$ induced by the underlying graph
functor $U:\Gpd\to\Gph$.
\end{proposition}

\begin{proof} By Proposition \ref{same connected components} it is sufficient to
prove that the inclusion functor $ \utor^2_U(\Pi,\tilde{A}_1)\hto
\Tor^2_U(\Pi,\tilde{A}_1) $ factors through the full and faithful functor of
Proposition \ref{2-exts of gpd are tor}. This follows from the fact that the
free groupoid on a graph has as objects the vertices of the graph and that the
counit map for the free adjunction is the identity on objects. This implies that
the fiber groupoid of any $U$-split 2-torsor in $\utor^2_U(\Pi,\tilde{A}_1)$ is
actually a 2-groupoid and then, by Proposition \ref{xm as groupoids}, it is
isomorphic to the groupoid associated by the functor \gpd\ to a crossed module
which is the fiber of a 2-extension of $\Pi$ by $A$.
\end{proof}

\gdef\alginv#1#2#3#4{#4 $\eta_{#2}:P_{#2}(\cbc)\to P_{#1}(\cbc)$ uniquely
corresponds to an element $$ k_{#2} \in H^2_{\cotrip[#1]} \bpar{P_#1(\cbc),
\tilde{A}_#1} $$ where $A = \pi_{#2} \bpar{P_#1(\cbc)}$. This cohomology element
will be called the algebraic #3 \pin\ of the crossed complex $\cbc$.}
\alginv12{second}{The last results show that the fibration}

\bigskip

For $n>1$, an observation about $q_n:\c_n\rightarrow\c_n/\Img\partial_{n+1}$
similar to the one made for $q_1$ holds, namely that $q_n$ is (not only the
cokernel of $\partial_{n+1}$, but also) the cokernel of
$\bar{\partial}_{n+1}:\c_{n+1}/\Img\partial_{n+2}\rightarrow\c_n$. As a
consequence, $\eta_{n+1}$ \bpar{described by diagram \eqref{eta n+1}} represents
a 2-extension of the $n$-crossed complex $P_n(\cbc)$ by the $\pi_1(\cbc)$-module
$\pi_{n+1}(\cbc)$, according to the following definition, which extends that of
2-extensions of groupoids:

\begin{definition}
\label{ext of crsn} If
$\cbc=(\c_n\xrightarrow{\partial_n}\cdots\xrightarrow{\partial_2} \termxm[\bg])$
is an $n$-crossed complex with $n>2$, and $A:\Pi\rightarrow\Ab$ is a
$\Pi$-module over the fundamental groupoid $\Pi=\pi_1(\cbc)$ of \cbc, a
2-extension of \cbc\ by $A$ is an exact sequence in the category of \bg-groups.
$$
\everyentry={\vphantom{\big(}} \xymatrix@C=1.25pc@R=1.5pc {0\ar[r] & A\ar[r] &
E_1 \ar[r]^\sigma & E_0\ar[r]^\tau & C_n\ar[r] & 0,}
$$ (where $A$ is considered as a \bg-module via the canonical projection
$q:\bg\rightarrow\Pi$), such that
\begin{equation}
\label{n+1crs} \xymatrix@C=1.25pc@R=1.5pc {\e_1\ar[r]^\sigma & \e_0
\ar[r]^-{\partial_n\tau} & \c_{n-1}\ar[r] &\cdots\ar[r]&\c_2\ar[r]^{\partial_2}
& \termxm[\bg]}
\end{equation} is an $(n+1)$ \bg-crossed complex, where $\e_i=(\bg,E_i,0),\;
i=0,1$. For $n=2$, we define a 2-extension of a crossed module
$\c=(\bg,C,\delta)$ by a $\Pi$-module $A$ as an exact sequence of \bg-groups
\begin{equation}
\label{2ext xm} \everyentry={\vphantom{\Big(}} \xymatrix@C=1.25pc@R=1.5pc
{0\ar[r] & A \ar[r] & E_1 \ar[r]^\sigma & E_0\ar[r]^\tau & C\ar[r] & 0,}
\end{equation} such that
\begin{equation}
\label{2ext xm2} \xymatrix@C=1.25pc@R=1.5pc {\e_1\ar[r]^\sigma &
\e_0\ar[r]^-{\delta\tau} & \termxm[\bg]}
\end{equation} is a 3-crossed complex, where $\e_0=(\bg,E_0,\delta\tau)$ and
$\e_1=(\bg,E_1,0)$.
\end{definition}

Let us note that, as in Definition \ref{2exgpd}, a 2-extension of a $n$-crossed
complex \cbc\ can be seen as an exact sequence
\begin{equation}
\everyentry={\vphantom{\Big(}} \xymatrix@C=1.25pc@R=1.5pc {0\ar[r] & \cba \ar[r]
& \cbe_1 \ar[r]^\sigma & \cbe_0\ar[r]^\tau & \cbc\ar[r] & 0,}
\end{equation} in the category $\bpar{\crs_n}_{T_{n-1}(\cbc)}$ of $n$-crossed
complexes with a fixed $(n-1)$-truncation, together with an extra structure in
the central part that makes \eqref{n+1crs} or \eqref{2ext xm2} an
$(n+1)$-crossed complex. The $n$-crossed complexes $\cba$ and $\cbe_i$, $i=0,1$,
have at dimension $n$ the crossed modules $\zero(A) = (\bg,A,0)$ and $\e_i$
respectively.

Using those definitions, our discussion can be summarized in the following
statement:
\begin{proposition} For all $n\geq 1$ the fibration $\eta_{n+1}$ of the
Postnikov tower of $\cbc$ provides a 2-extension of $P_n(\cbc)$ by the
$\pi_1(\cbc)$-module $\pi_{n+1}(\cbc)$.
\end{proposition}

\begin{proof} For $n=1$, the 2-extension associated to $\eta_2$ is given by the
crossed module $P_2(\cbc)$ and the sequence \eqref{k2 ex seq gpd}. For $n>1$ the
2-extension associated to $\eta_{n+1}$ is given by the sequence
$$
\everyentry={\vphantom{\big(}} \xymatrix@C=1.25pc@R=1.5pc {0\ar[r] &
\pi_{n+1}(\cbc)\ar[r] & C_{n+1}/\Img\partial_{n+2}
\ar[r]^-{\bar{\partial}_{n+1}} & C_n \ar[r]^-{q_n} &
C_n/\Img\partial_{n+1}\ar[r] & 0}.
$$

\end{proof}

As suggested above, it is not difficult to extend this proposition to the case
$n=0$ by giving an appropriate definition of a 2-extension of a set $X$ by an
$X$-indexed family of groups. The details are in \cite{Garcia2003}.

As  in the case of groupoids (case $n=1$), the 2-extensions of an $n$-crossed
complex \cbc\ by a fixed $\pi_1(\cbc)$-module $A$ constitute a category, denoted
$\Ext^2(\cbc,A)$, where morphisms between extensions are defined in the obvious
way, that is, as commutative diagrams
$$
\everyentry={\vphantom{\big(}} \xymatrix@C=1.25pc@R=1.5pc {0\ar[r] &
A\ar@{=}[d]\ar[r] & E_1\ar[d] \ar[r]^\sigma & E_0\ar[d]\ar[r]^\tau &
C_n\ar[r]\ar@{=}[d] & 0, \\ 0\ar[r] & A \ar[r] & E_1'
\ar[r]^{\sigma'} & E_0'\ar[r]^{\tau'} & C_n\ar[r] & 0.}
$$ Again, $\Ext^2[\cbc,A]$ denotes the set of connected components of
$\Ext^2(\cbc,A)$.

\bigskip

It is now necessary a more elaborate analysis than the one made for groupoids in
order to show that the 2-extensions of crossed modules  defined above represent
cohomology elements in a cotriple cohomology of crossed modules (of course the
one corresponding to the cotriple induced by $U_2$, Proposition \ref{pxm
tripleable over agpd}).

Let $\c=(\bg,C,\delta)\in\crs_2=\xm$, let $\Pi=\pi_1(\c)$ be its fundamental
groupoid and let $A:\Pi\rightarrow\Ab$ be a system of local coefficients. It is
not difficult to prove that the resulting crossed module $\abgc A2$ has base
groupoid equal to $\bg$ and structure $\bg$-group the functor (also denoted
$\abgc A2$ by abuse of notation) $\abgc A2:\bg\to \Gr$ given by the cartesian
product of groups, $\abgc A2(x) = C(x)\times A(x)$, with action
$\laction{u}{(v,b)} = (\laction uv,\laction ub)$. Furthermore, the crossed
module connecting morphism of $\abgc A2$, $\delta'$, has components
$\delta'_x:C(x)\times A(x)\to\End_\bg(x)$ given by $\delta'_x\bpar{(u,a)} =
\delta(u)$. Regarding $\abgc A2$ as an internal abelian group object in
$\xm/\c$, it will be taken as a system of global coefficients.\label{ab grp in
xm}

Let us consider now a 2-extension of \c\ by $A$ such as \eqref{2ext xm}. By the
condition that diagram \eqref{2ext xm2} be a 3-crossed complex, the action of
$\Img \delta\tau$ on $\widehat{E}_1$ is trivial and the cartesian product
$\e_0\times\e_1$ in $\xm_{\bg}$ is given by
$$
\e_0\times\e_1 =(\bg,E_0\times E_1,\delta\tau p_0).
$$ This crossed module has two obvious maps of crossed modules to $\e_0$,
namely, the canonical projection $p_0:\e_0\times\e_1\rightarrow \e_0$ and the
map determined by $t:(x,y)\mapsto x\sigma(y)$ from $E_0\times E_1$ to $E_0$.
These two maps have a common section determined by the map $x\mapsto (x,0)$ and
the resulting internal graph in \xm,
$$
\cbe: \xymatrix@C=1pc@R=1.75pc
{**[l]{\e_0\times\e_1}\ar@<.7ex>[rr]|-{\phantom{.}p_0\phantom{.}}\ar@<-0.5ex>[rr]_-t
& & \e_0 \ar@/_/@<-1.2ex>[ll] }
$$ admits a unique groupoid structure in which the multiplication is determined
by
$$
\langle(x,y),(x',y')\rangle\mapsto (x,yy').
$$

Since the crossed module of connected components of \cbe\ (the coequalizer of
$p_0$ and $t$) is easily verified to be the canonical map $(\tau,1_\bg):\e_0\to
\c$ determined by $\tau:E_0\rightarrow C$ (and the identity of \bg\ as change of
base), we can take the structure of internal groupoid of \cbe\ as the fiber
groupoid of a $(\abgc A2,2)$-torsor above \c. To define such a torsor we just
need to give the corresponding cocycle map $\alpha:\End(\cbe)\rightarrow \abgc
A2$. The domain of this map is the crossed module obtained as the equalizer of
$p_0$ and $t$. Since the condition defining this equalizer is $x=x\sigma(y)$,
one quickly finds that $\End(\cbe)=(\bg,E_0\times\ker \sigma,\delta\tau p_0)$.
The required morphism $\alpha$ is defined as the map of \bg-crossed modules
determined by the following map of \bg-groups:
$$ E_0\times\ker\sigma\xrightarrow{\tau\times\bar{\sigma}} C\times (A\circ q),
$$ where $\bar{\sigma}$ is the canonical isomorphism $\ker\sigma\cong A\circ q$
induced by the exactness of $0\rightarrow A\circ q\rightarrow
E_1\xrightarrow{\sigma} E_0$.

The above arguments have prepared the ground for the following:

\begin{proposition}
\label{2-ext of xm are tor} For any crossed module $\c=(\bg,C,\delta)$ and any
$\pi_1(\c)$-module $A:\pi_1(\c)\rightarrow\Ab$, there is a full and faithful
functor
$$
\Ext^2(\c,A)\rightarrow \Tor^2_{U_2}(\c,\abgc A2),
$$ where $U_2$ is the monadic functor of Proposition \ref{pxm tripleable over
agpd}.
\end{proposition}

\begin{proof} It is a simple exercise to verify that the construction given
above indeed produces a 2-torsor above $\c$ with coefficients in $\abgc A2$ from
any 2-extension of $\c$ by $A$, and that this construction is functorial.
Furthermore, if one examines the correspondence between morphisms of extensions
$$
\vcenter{\everyentry={\vphantom{\big(}} \xymatrix@C=1.25pc@R=1.5pc {0\ar[r] &
A\ar@{=}[d]\ar[r]^j & E_1\ar[d]_{f_1} \ar[r]^\sigma & E_0\ar[d]_{f_0}\ar[r]^\tau
& C \ar[r]\ar@{=}[d] & 0, \\ 0\ar[r] & A \ar[r]^{j'} & E_1'
\ar[r]^{\sigma'} & E_0'\ar[r]^{\tau'} & C\ar[r] & 0 }} \qquad
\begin{aligned}
\tau' \f_0&= \tau \\ f_0\sigma &=\sigma' f_1 \\ f_1 j &=j'
\end{aligned}
$$
and morphisms of 2-torsors
\vglue -7 mm
$$
\everyentry={\vphantom{\big(}} \xymatrix{& & E_0 \times \ker \sigma
\ar[dd]_{f_0\times f_1} \ar[dl]_{p_0} \ar[ddr]^{\bar{\alpha}} & \\
**[l]E_0\times E_1 \ar@<.7ex>[r]|-{\phantom{.}p_0}\ar@<-0.7ex>[r]_-t
\ar@<-3ex>[dd]_{f_0\times f_1} & E_0 \ar@/_1pc/[l] \ar[dd]_{f_0}
\ar[ddr]|<>(.4){\phantom{-}\sigma\phantom{-}} & \\ & & E_0'\times \ker \sigma'
\ar[dl]|-{\phantom{.}p_0'} \ar[r]^-{\bar{\alpha}'} & **[r]C\times (A\circ q)
\ar[dl]^{p_C} & \\ **[l]E_0'\times E_1'\labelmargin{2ex}
\ar@<.7ex>[r]|-{\phantom{.}p_0'}\labelmargin{.5ex} \ar@<-0.7ex>[r]_-{t'} & E_0'
\ar[r]_{\sigma'}\ar@/_1pc/[l] & C }
\mskip - 120 mu
\begin{aligned}
\\[7 mm]
(1_{\tilde{A}}\times f_0) \alpha &=\alpha'  (f_0\times f_1) \\
f_0 t &=t' (f_0\times f_1)
\end{aligned}
$$
it becomes evident that it is a bijective correspondence and therefore
the functor from 2-extensions of \c\ by $A$ to $(\abgc A2,2)$-torsors
above \c\ is full and faithful.

It only remains to prove that the internal groupoid \cbe\ in \xm\ that we have
associated to an extension of \c\ by $A$ is $U_2$-split. For that, let us
consider $U_2(\tau,1_\bg)=(\tau,1_\bg)$, where $\tau$ in the right-hand side is
regarded a map of sets from the set of arrows of $\widehat{E}_0$ to the set of
arrows of $\widehat{C}$. This map of sets is surjective and for any section $s$
of this map, the pair $(s,1_\bg)$ is a section of $(\tau, 1_\bg)$ in \agpd.

We need to consider also the map $\langle p_0,t\rangle :\e_0\times\e_1\rightarrow\e_0\times_\c
\e_0$ and prove that $U_2(\langle p_0,t\rangle)$ has a section. Giving a section for this map is
equivalent to giving, for each object $x\in\bg$, a map of sets $K_x:E_0(x)\times_{C(x)}
E_0(x)\rightarrow E_1(x)$ such that $\sigma_x(K_x(u,v))=u^{-1}v$. Taking into account that
$E_0(x)\times_{C(x)}E_0(x)$ is the set of pairs $(u,v),\; u,v\in E_0(x)$ such that
$\tau_x(u)=\tau_x(v)$, it is clear that for every $(u,v)\in E_0(x)\times_{C(x)} E_0(x)$ we have
$u^{-1}v\in\ker \tau_x=\Img \sigma_x$. We can choose any section $\beta$ of the set map
$E_1(x)\xrightarrow{\sigma_x}\Img\sigma_x$ and define $K_x(u,v)=\beta(u^{-1}v)$. This proves
$U_2(\langle p_0,t\rangle)$ has a section, and therefore \cbe\ is $U_2$-split.
\end{proof}

\begin{proposition} For any crossed module $\c$ and any $\pi_1(\c)$-module
$A:\pi_1(\c) \to \Ab$,
$$ H_{\cotrip[2]}^2(\c,\abgc A2)\cong \Ext^2[\c,A].
$$ where $\abgc A2$ is the abelian group object in $\xm$ obtained from $A$ as
indicated in page \pageref{ab grp in xm}.
\end{proposition}

\begin{proof} This is consequence of Proposition \ref{same connected components}
and the fact that the inclusion $\utor_{U_2}^2(\c,\abgc
A2)\hto\Tor_{U_2}^2(\c,\abgc A2)$ factors through the full and faithful
inclusion of Proposition \ref{2-ext of xm are tor}. To see the last statement,
observe that the counit map for the free adjunction $F_2:\agpd\leftrightarrows
\xm:U_2$ is always the identity at the level of base groupoid. Therefore the
fiber groupoid (internal in \xm) of any 2-torsor in $\utor_{U_2}^2(\c,\abgc A2)$
lives in $\xm_\bg$. Finally it is straightforward to see that any internal
groupoid in $\xm_\bg$ is isomorphic to the groupoid built from a 2-extension as
in the proof of Proposition \ref{2-ext of crsn are tor}.
\end{proof}

\alginv23{third}{The last results show that the fibration}

\subsubsection*{The higher invariants.}

Our next objective is to establish a general bijection between
$\Ext^2[P_n(\cbc),\pi_{n+1}(\cbc)]$ and the set of elements in certain cotriple cohomology in the
category of $n$-crossed complexes for $n\geq3$. An essential step in this process is determining
the coefficients to be used to calculate the cohomology. These coefficients, (constituting an
internal abelian group object in the category $\crs_n/P_n(\cbc)$) are obtained from the
``homotopy group'' $\pi_{n+1}(\cbc)$ using the fact that $\pi_{n+1}(\cbc)$ is a module over the
fundamental groupoid of $P_n(\cbc)$ (equal to the fundamental groupoid of $\cbc$, in turn equal
to $P_1(\cbc)$). In general, if $\cbc\in\crs_n$ is a $n$-crossed complex for $n>1$,  pulling back
along the canonical map $\cbc \to \Pi=\pi_1(\cbc) = P_1(\cbc)$ (a finite product-preserving
functor) produces an abelian group object in $\crs_n/\cbc$ from any abelian group object in
$\crs_n/\Pi$. This allows us to reduce the search for our coefficients to obtaining an abelian
group object in $\crs_n/\Pi$. \guardar{This can be done in several ways. For our purposes what we
need is the abelian group object defined in the following way. For $n=1$ we take the canonical
abelian group structure on the semidirect product $(\semi\Pi A\leftrightarrows\Pi)\in\Gpd/\Pi =
\crs_1/\Pi$. If $n>1$ we use} For every $n\geq2$ there is a functor $\ins_n:\Ab^\Pi\to\crs_n$
taking a $\Pi$-module $A:\Pi\to\Ab$ to the $n$-crossed complex over $\Pi$
\begin{equation}
\label{insert n} \ins_n(A) = (\zero(A)\to\initxm[\Pi]\to
\cdots\to\initxm[\Pi]\to\termxm[\Pi]).
\end{equation} This functor can be regarded as taking its values in $\crs_n/\Pi$
via the obvious map $(1_\Pi,{\boldsymbol 0}):\ins_n(A)\to\Pi$, and when regarded
this way it becomes a finite product preserving functor whose value at a
$\Pi$-module $A$ will be denoted $\abgp An$. Since (the theory of abelian groups
being commutative) $A$ has the structure of an internal abelian group object in
$\Ab^\Pi$, the fact that $\ins_n:\Ab^\Pi\to\crs_n/\Pi$ preserves finite products
implies that $\abgp An$ has a structure of internal abelian group object in
$\crs_n/\Pi$.

If $\cbc$ is fixed by the context and $\Pi = \pi_1(\cbc)$, the abelian group
object in $\crs_n/\cbc$ obtained from a $\Pi$-module $A$ after pulling
$\ins_n(A)\to\Pi$ back along the canonical map $\cbc\to\Pi$ will be denoted
$\abgc An$, so that we have,
\begin{equation}
\label{abelian group in xm} \everyentry={\vphantom{\big(}} \vcenter{
\xymatrix@C=1.25pc@R=1.5pc {\abgc An\ar[d]\ar@{}[dr]|{\text{pb}} \ar[r] &
\ins_n(A) \ar[d] \\ \cbc\ar[r]_-{\text{can.}} & \Pi} }
\end{equation}

Let $n \geq 3$, \cbc\ an $n$-crossed complex, $\Pi=\pi_1(\cbc)$ and
$A:\Pi\rightarrow\Ab$ a $\Pi$-module. We define $\abgc An$ by \eqref{abelian
group in xm} and take the resulting abelian group object in $\crs_n/\cbc$ (also
denoted $\abgc An$) as a system of global coefficients for 2-torsors.

\begin{proposition}
\label{2-ext of crsn are tor} Let $n \geq 3$, for any $n$-crossed complex \cbc\
and any $\pi_1(\cbc)$-module $A$ there is a full and faithful functor
$$
\Ext^2(\cbc,A)\rightarrow \Tor^2(\cbc,\abgc An).
$$
\end{proposition}

\begin{proof} Let us consider a 2-extension of \cbc\ by $A$,
$$
\everyentry={\vphantom{\big(}} \xymatrix@C=1.25pc@R=1.75pc {0\ar[r] & A\ar[r] &
E_1 \ar[r]^\sigma & E_0\ar[r]^\tau & C_n\ar[r] & 0,}
$$ let us also denote $\e_i=(\bg,E_i,0)$ as in Definition \ref{ext of crsn}. The
\bg-module $E_0\oplus E_1$ with the zero map gives a \bg-crossed module
$\e_0\oplus\e_1$ which, substituted for $\c_n$ in \cbc\ with the boundary map
$\partial_n\tau p_0:\e_0\oplus \e_1\rightarrow \c_{n-1}$, gives rise to an
$n$-crossed complex over \bg,
$$
\cbf_1: \e_0\oplus \e_1\xrightarrow{\partial_n\tau p_0}
\c_{n-1}\longrightarrow\cdots\longrightarrow\c_2\longrightarrow \termxm[\bg].
$$ We will take this as the object of arrows of an internal groupoid in
$\crs_n$. As the object of objects we take the $n$-crossed complex
$$
\cbf_0: \e_0\xrightarrow{\partial_n\tau}
\c_{n-1}\longrightarrow\cdots\longrightarrow\c_2\longrightarrow \termxm[\bg].
$$ The ``source" map $s:\cbf_1\rightarrow\cbf_0$ is the obvious map of crossed
complexes induced by the projection $p_0:E_0\oplus E_1\rightarrow E_0$, and the
``target" map $t:\cbf_1\rightarrow\cbf_0$ is the one induced by the map $x\oplus
y\mapsto x\sigma(y)$ form $E_0\oplus E_1$ to $E_0$. Then, the canonical
inclusion $E_0\hookrightarrow E_0\oplus E_1$ determines a common section for $s$
and $t$, and we obtain an internal groupoid in $\crs_n$ in which composition is
determined by the map $(x\oplus y,x'\oplus y')\mapsto x\oplus
yy'.$

It is a simple matter to show that the $n$-crossed complex of endomorphisms of
this groupoid (the equalizer of $s$ and $t$ in $\crs_n$) is
$$
\cbe: E_0\oplus (A\circ q) \longrightarrow
\initxm[\bg]\longrightarrow\cdots\longrightarrow \initxm[\bg]\longrightarrow
\termxm[\bg],
$$ and the $n$-crossed complex of connected components of this groupoid (the
coequalizer of $s$ and $t$ in $\crs_n$) is \cbc. Thus, the above internal
groupoid could be taken as the fiber of a torsor in $\Tor^2(\cbc,\abgc An)$ if a
cocycle map $\alpha:\cbe\rightarrow\cbc\times\abgc An$ can be given. A simple
calculation shows that $\cbc\times\abgc An$ is the $n$-crossed complex
$$
\left( C_n\oplus (A\circ q), 0\right)\longrightarrow
\c_{n-1}\longrightarrow\cdots\longrightarrow\c_2\longrightarrow \termxm[\bg]
$$ and the cocycle map $\alpha$ can be defined as the obvious map induced by
$$
\tau\oplus (A\circ q):E_0\oplus (A\circ q)\longrightarrow C_n\oplus (A\circ q),
$$ note that the square
$$
\xymatrix@C=10ex{ E_0\oplus (A\circ q) \ar[d] \ar@{}[dr]|{\text{pb}}
\ar[r]^-{\tau\oplus (A\circ q)} & C_n\oplus (A\circ q) \ar[d] \\ E_0
\ar[r]_-\tau & C_n }
$$ is a pullback in the category of \bg-modules.

This defines the functor $\Ext^2(\cbc,A)\rightarrow \Tor^2(\cbc,\abgc An)$ on
objects. On morphisms the functor is defined by the obvious map, which
establishes a bijection between the morphisms between two extensions and the
morphisms of torsors between the corresponding torsors. In this way one gets
the desired full and faithful functor.
\end{proof}

We will define now a monadic functor $U_n$ in $\crs_n$ and prove that the torsors obtained by the
functor of Proposition \ref{2-ext of crsn are tor} are $U_n$-split. Let $\acrs_{n-1}$ be the
category whose objects are triples $(X,f,\cbc)$ where $X$ is a set, \cbc\ is an ${(n-1)}$-crossed
complex, and $f:X\rightarrow \arr (\widehat{\pi_{n-1}(\cbc)})$ is a map from $X$ to the set of
arrows of the totally disconnected groupoid associated to $\pi_{n-1}(\cbc):\pi_1(\cbc)\rightarrow
\Ab$. An arrow $(X,f,\cbc)\rightarrow(X',f',\cbc')$ in $\acrs_{n-1}$ is a pair
$(\alpha,\beta)$ where $\alpha:X\rightarrow X'$ is a map of sets and
$\beta:\cbc\rightarrow\cbc'$ is a map in $\crs_{n-1}$ such that $\tilde{\beta} f = f'
\alpha$ (where $\tilde{\beta}:\pi_{n-1}(\cbc)\rightarrow\pi_{n-1}(\cbc')$ is the obvious
natural map induced by $\beta$).

There is an obvious forgetful functor $U_{n}:\crs_{n}\rightarrow\acrs_{n-1}$,
taking an $n$-crossed module \cbc\ to the triple
$\big(\arr(\widehat{C}_{n}),\widehat{\partial}_{n}, T_{n-1}(\cbc)\big)$ (the
``simple truncation" functor $T_{n-1}$ is defined in section \ref{crs}).

\begin{proposition}
\label{crsn tripleable over acrs} For $n> 2$, the forgetful functor
$U_{n}:\crs_{n}\rightarrow\acrs_{n-1}$ is monadic.
\end{proposition}

\begin{proof} We begin by defining a left adjoint. This is very similar to the
case of Proposition \ref{pxm tripleable over agpd}. The main difference is that in this case we
directly get a crossed module with trivial connecting morphism instead of a simple pre-crossed
module. Given $(X,f,\cbc)$ in $\acrs_{n-1}$, with
$$
\cbc: \c_{n-1}\xrightarrow{\partial_{n-1}}\c_{n-2}
\longrightarrow\cdots\longrightarrow\c_2\xrightarrow{\partial_2} \termxm[\bg]
$$ and $\c_2=(\bg,C,\delta)$, we need to define a trivial \bg-crossed module
$\c_{n}=(\bg,C_{n},0)$ on which $\Img(\delta)$ acts trivially. For this, it is
sufficient to define a $\pi_1(\cbc)$-module
$\overline{C}_{n}:\pi_1(\cbc)\rightarrow\Ab$ and to put
$C_{n}=\overline{C}_{n}\circ q$, where $q:\bg\rightarrow\pi_1(\cbc)$ is the
canonical map. For each object $x\in\bg$ we define the abelian group
$$
\overline{C}_{n}(x)= {\sff}_{\text{ab}}\bigg( \coprod_{z\in\obj(\bg)}
\Big(\pi_1(\cbc)(z,x)\times\coprod_{v\in\pi_{n-1}(\cbc)(x)} \fbr(f,v)\Big)
\bigg).
$$ that is, $\overline{C}_{n}(x)$ is the free abelian group generated by all
pairs $\langle t,u\rangle$ where $t:z\to x$ is a map in $\pi_1(\cbc)$ and $u\in
X$ such that $f(u)\in \pi_{n-1}(\cbc)(z)$.

Given $s:x\to y$ in $\pi_1(\cbc)$, $\overline{C}_{n}(s)$ is the homomorphism
defined on the generators of $\overline{C}_{n}(x)$ by
$$
\overline{C}_{n}(s)\big({\langle t,u\rangle}\big) = {\langle st,u\rangle}.
$$ This defines the functor $\overline{C}_{n}$. We now take $C_{n} =
\overline{C}_{n}\circ q$ and we obtain a \bg-crossed module $\c_{n} = (C_{n},0)$
on which (by construction) $\Img(\delta)$ acts trivially. Thus, to obtain an
object in $\crs_{n}$ we only need a morphism $\partial_{n}:\c_{n}\to\c_{n-1}$
such that $\partial_{n-1}\partial_{n} = 0$. Given $x\in\bg$ we define
$(\partial_{n})_x:C_{n}(x)\to C_{n-1}(x)$ as the group homomorphism defined on
the generators of $C_{n}(x)$ by
$$ (\partial_{n})_x({\langle t,u\rangle}) = \laction t{f(u)}
=C_{n-1}(t)\big(f(u)\big).
$$

We have thus defined an $n$-crossed complex
$$ F_{n}\big((\cbc, X, f)\big):\quad\c_{n} \xto{\partial_{n}} \c_{n-1}
\xto{\partial_{n-1}} \c_{n-2}\to \cdots \to \c_2 \xto{\partial_2}
\c_1=\termxm[\bg],
$$ It remains to define $F_{n}$ on arrows. Given an arrow
$(\alpha,\beta):(X,f,\cbc)\to(X',f',\cbc')$ in $\acrs_{n-1}$, a map of crossed
complexes is determined by $\alpha$ with the additional component $\alpha_{n}$
being the natural transformation with components defined on the generators of
$\c_{n}(x)$ by
$$
\big(\alpha_{n}\big)_x(\langle t,u\rangle) = \langle\alpha_0(t),\beta(u)\rangle,
$$ where $\alpha_0$ is the ``change of base groupoid'' part of $\alpha$. It is
now a straightforward exercise to verify that $F_{n}$ is a functor
$\crs_{n}\to\acrs_{n-1}$ and that this functor is left adjoint to the forgetful
$U_{n}$ (see \cite{Garcia2003}{, Proposici\'on 3.2.10, pp. 161--166} for the
details).

As in Proposition \ref{pxm tripleable over agpd} is is easy to see that $U_n$
reflects isomorphisms. The proof is completed by a calculation similar to the
one used in the said proposition, showing that $U_n$ preserves coequalizers of
$U_n$-contractible pairs.
\end{proof}

Note that the counit $\epsilon_\cbc:F_{n}U_{n}(\cbc) \to \cbc$ of the adjunction
$F_{n}\dashv U_{n}$ is an identity at dimensions other than $n$. The same thing
occurs with the unit $\eta_\cbc : U_{n}F_{n}U_{n}(\cbc) \to U_{n}(\cbc)$ and
with its image by $F_{n}$.

\begin{proposition}
\label{2-ext of crsn are uspl tor} In the hypothesis of Proposition \ref{2-ext
of crsn are tor}, the 2-torsor obtained from any 2-extension of $\cbc$ by $A$ is
$U_n$-split. Therefore, the functor defined in Proposition \ref{2-ext of crsn
are tor} actually represents a full and faithful functor
$$
\Ext^2(\cbc,A)\rightarrow \Tor_{U_n}^2(\cbc,\tilde{A}).
$$
\end{proposition}

\begin{proof} Since all components of the coequalizer map $\cbf_0 \to \cbc$ in
dimensions $<n$ are identities, that the obtained torsor is $U_n$-split is
equivalent to the fact that the $n$-dimensional component $\tau: F \to C_n$ is
surjective.
\end{proof}

We denote $\cotrip[n]$ the cotriple induced on $\crs_n$ by the the monadic
functor $U_n$, that is, $\cotrip[n] = F_nU_n$. Given a $n$-crossed complex \cbc,
the adjoint pair $(F_n,U_n)$ induces an adjoint pair (also denoted  $(F_n,U_n)$)
between the corresponding slice categories $\crs_{n}/\cbc$, and
$\acrs_{n-1}/U_n(\cbc)$. Furthermore, the induced $U_n$ on the slice is again
monadic. Thus we obtain again a cotriple on the category $\crs_{n}/\cbc$ and
this cotriple will also be denoted $\cotrip[n]$.

\begin{proposition} For any $\cbc\in\crs_n$ and any $A:\pi_1(\cbc) \to \Ab$, the
inclusion functor $\utor^2_{U_n}(\cbc,\abgc An) \hto {\Tor}^2_{U_n}(\cbc,\abgc
An)$ factors through the full and faithful functor of Proposition \ref{2-ext of
crsn are uspl tor}. As a consequence,
$$ H_{\cotrip[n]}^2(\cbc,\abgc An)\cong \Ext^2[\cbc,A].
$$
\end{proposition}

\begin{proof} Again this proof is based on the fact that the counit of the
adjunction $F_n\dashv U_n$ is the identity on the $(n-1)$ truncation, then the
groupoid fiber of any 2-torsor in $\utor^2_{U_n}(\cbc,\abgc An)$ is isomorphic
to one which comes from a 2-extension.
\end{proof}

In view of the last results, the fibrations in the \pto\ of a \cc\ \cbc\ can be
regarded as 2-extensions, and we have seen they represent a cotriple cohomology
of crossed complexes. \alginv n{n+1}{$(n+1)^\text{th}$}{The fibration}

\subsection{Singular cohomology and the topological invariants}

\label{sec the sing coh}

The knowledge of the above algebraic Postnikov invariants is sufficient
information to reconstruct the (homotopy type of the) crossed complex \cbc. This
solves the problem of calculating the Postnikov invariants of any space having
the homotopy type of a crossed complex. However, in order to relate the above
invariants to the usual calculation of the Postnikov invariants of a space (as
elements in the singular cohomology of the space), we would like to show that
our algebraically obtained invariants determine, in a natural way, singular
cohomology elements.

Remember that the  ``fundamental crossed complex" functor (see \cite{BrHi1991})
associates a crossed complex to each simplicial set,
$$
\Pi:\sset\rightarrow\crs.
$$ This functor has a right adjoint ``nerve of a crossed complex" which
associates to each crossed complex \cbc\ the simplicial set whose set of
$n$-simplices is the set of maps of crossed complexes
$$
\ner(\cbc)_n=\crs\bpar{\Pi(\Delta[n]),\cbc}.
$$ Applying this nerve functor to a fibration of crossed complexes one obtains a
fibration of simplicial sets whose fibers have the same homotopy type as those of the initial
fibration. \vskip0pt Since the geometric realization functor $\sset\rightarrow\Top$ also
preserves fibrations and the homotopy type of their fibers, we can define a functor ``classifying
space of a crossed complex" which again has the same property. If we restrict ourselves to the
category \t\ of those spaces having the homotopy type of a crossed complex we obtain
two functors \begin{equation}\label{fund crs ladj to geom realiz}
\Pi:\t \to \crs \,,\qquad B: \crs \to \t
\end{equation} which induce an equivalence in the corresponding homotopy
categories and which allow us to reduce the calculation of the Postnikov towers
of the spaces in \t, to the calculation of Postnikov towers in \crs.

We define the singular cohomology of a crossed complex \cbc\ as the singular
cohomology of its geometric realization, $B(\cbc)$ or, equivalently, of the
geometric realization of its nerve. Thus, if $\Pi = \pi_1(\cbc)$ and
$A:\Pi\to\Ab$, we define
$$
\scoH^m(\cbc,A) = \scoH^m\bpar{B\bpar{\ner(\cbc)},A},
$$ and our aim is to establish a natural map
$$ H^2_{\cotrip[n]}\bpar{P_n(\cbc),\abgc An} \xto{\ \
}\scoH^{n+2}\bpar{P_n(\cbc), A}.
$$ In order to establish this map we will obtain a representation of the
singular cohomology of an $n$-crossed complex by homotopy classes of maps in a
certain category of simplicial $n$-crossed complexes (denoted $\scrs{n}$) which
is naturally associated to the cotriple $\cotrip[n]$. Then the desired map will
be induced by a natural map from Duskin's representation of the cotriple
cohomology to our representation in $\scrs{n}$ of the singular cohomology. This
representation is actually obtained as a ``lifting'' of the generalized
Eilenberg-MacLane representation in simplicial sets (see \cite{GoJa1999} and
\cite{BuFaGa2002}):
\begin{equation}\label{gener eilmacl}
\scoH^m\bpar{\cbc,A} =
\big[\ner(\cbc), L_\Pi(A,m)\big]_{\sset/\ner(\Pi)}.
\end{equation} where $L_\Pi(A,m)$ is the canonical fibration from the homotopy
colimit of the functor $K\bpar{A(\cdot),m}$ to $K(\Pi,1)$ ($= \ner(\Pi)$) (see
\cite{GoJa1999} and \cite{BuFaGa2002}). The ``ladder'' through which this
lifting is achieved is a chain of adjunctions
$$
\xymatrix{\sset\ar@<.5 ex>[rr]^-{\f_1=G} && \scrs{1} \ar@<.5
ex>[ll]^-{\W_1=\ner} \ar@<.5 ex>[r]^-{\f_2} &\scrs 2 \ar@<.5 ex>[l]^-{\W_2}
\quad\cdots\quad  \scrs{n-1}\ar@<.5 ex>[r]^-{\f_n} &\scrs{n} \ar@<.5
ex>[l]^-{\W_n}\ \cdots }
$$ having sufficiently good properties so as to preserve homotopy classes of
certain maps.

We will first establish the adjunctions and their properties in the cases $n=1$ (essentially done
in \cite{DwKa1984}) and $n=2$, which are the harder cases. Later, the cases $n\geq 3$ will be
dealt with.

Let us begin by introducing the categories $\scrs{n}$, $n\geq0$, whose
definition is motivated by a special property of the simplicial \ccs\ arising
from the cotriple $\cotrip[n]$.

Let $\cotrsr : \crs_n \to \dsimpl{{\crs_n}}$ denote the functor associating to
each \ncc\ \cbc\ the simplicial resolution of \cotrip[n] in \cbc, so that
$\cotrsr(\cbc)$ is a simplicial \cc\ which at dimension $k$ is equal to
$\cotrsr(\cbc)[k] = \cotrip[n]^{k+1}(\cbc)$ (where $\cotrip[n]^k$ is the
$k$-fold composition $\cotrip[n]\circ \cdots \circ\cotrip[n]$). We can state the
following (obvious) lemma, which provides the definition of $\scrs{n}$:

\begin{lemma}
\label{cotriple factors} For all $n\geq1$ and all $k\geq 0$ the \ccs\ $\cotrip[n]^k(\cbc)$ have
the same $(n-1)$-truncation and all faces and degeneracies of $\cotrsr(\cbc)$ at all dimensions
are morphisms of \ccs\ whose $(n-1)$-truncation is the identity of $T_{n-1}(\cbc)$. The functor
$\cotrsr$ factors through the subcategory $\scrs{n}$ defined by the following pullback of
categories:
\begin{equation}
\label{k2 ex seq xmm}
\vcenter{
\everyentry={\vphantom{\big(}}
\xymatrix@C=1.5pc@R=1.5pc { &\scrs{n}\ar[r] \ar@{^{(}->}[d]
\ar@{}[dr]|{\text{\rm pb}} &
\crs_{n-1}\ar@{^{(}->}[d]^{\text{\rm diag}}\\
{\crs_{n}}\ar[r]_-{\cotrsr}\ar@{.>}[ur]&{\dsimpl{{\crs_{n}}}}\ar[r]_-{(T_{n-1})_*}
& {\dsimpl{{\crs_{n-1}}}}} }
\end{equation}
\end{lemma}

The category $\scrs{1}$ is the category of those simplicial groupoids having the same set of
objects in all dimensions and all faces and degeneracies equal to the identity at the level of
objects. This is precisely the category (denoted \emph{Gd} in \cite{DwKa1984}) on which Dwyer and
Kan define the \emph{classifying complex} functor, $\W$, extending the classical functor
$\W:\dsimpl{\Gp} \to \dsimpl{\sets}$ of  \cite{EiMLa1954a} and \cite{May1967}. We will find it
convenient to define $\scrs{0} = \dsimpl{\sets}(=\sset)$ and to refer to $\W$ as $\W_1:\scrs{1}
\to \scrs{0}$. Note that $\scrs{1}$ can also be described as the category of groupoids
\emph{enriched} in simplicial sets.

The category $\scrs{2}$ is the subcategory of $\dsimpl{\xm}$ whose objects are simplicial crossed
modules $\Sigma:\op{\Delta} \to \xm$ such that for all dimensions the crossed modules $\Sigma_n$
have the same base groupoid and whose morphisms from $\Sigma$ to $\Sigma'$ are simplicial maps
$\alpha:\Sigma\to\Sigma'$ having in each dimension the same change-of-base functor.

\subsubsection{The lifting to $Gd = \scrs 1$.}

As we said above $\W_1$ is the functor $\W$ defined in \cite{DwKa1984}. We repeat here the
definition in order to have it handy. If $\Sigma\in\scrs{1}$ then the vertices of $\W_1(\Sigma)$
are the objects of $\Sigma$, the $n$-simplices ($n>0$) are the sequences of arrows
$(z_n\xto{u_{n-1}}z_{n-1}\to \cdots \to z_1\xto{u_0} z_0)$ where $u_i$ is an arrow in the
groupoid $\Sigma_i$ for $i=0,\dots,n-1$, and faces and degeneracies given by the formulas
\begin{align*} d_i(u_0,\dots,u_{n-1}) &= (u_0,\dots,u_{n-i-2},u_{n-i-1}\cdot
d_0u_{n-i},d_1u_{n-i+1},\dots,d_{i-1}u_{n-1}),\\ s_i(u_0,\dots,u_{n-1})
&=(u_0,\dots,u_{n-i-1},id, s_0u_{n-i},\dots,s_{i-1}u_{n-1}).
\end{align*} (Note the discrepancy with \cite[p.383]{DwKa1984} where there is an
error in the indices.)

Evidently, if $\Pi$ is a groupoid regarded as a constant simplicial groupoid, then the
$n$-simplices of $\W(\Pi)$ are precisely the $n$-simplices of $\ner (\Pi)$ and in fact we have
$\W_1(\Pi) = \ner(\Pi)$. Thus, $\W_1$ can be regarded as a functor
\begin{equation}\label{wbar}
\W_1:\scrs{1}/\Pi \longrightarrow \dsimpl{\sets}/\ner(\Pi).
\end{equation}

\begin{proposition} Let $\Pi$ be a groupoid, let $A:\Pi\to\Ab$ be a $\Pi$-module and
let $\textcolor{black}{\tilde{A}_1 = \ins_1(A) }=\semi\Pi A$, seen as an abelian
group object in
$\Gpd/\Pi$. Then the simplicial object
$K(\textcolor{black}{\tilde{A}_1},n)\in\dsimpl{\bpar{\Gpd/\Pi}}$ $= \dsimpl{\Gpd}/\Pi$ actually
belongs to $\scrs{1}/\Pi$ (where $\Pi$ is regarded as the trivial ``constant'' simplicial
groupoid) and \begin{equation}\label{w bar uno} \W_1\bpar{K(\textcolor{black}{\tilde{A}_1},n)} =
L_\Pi(A,n+1).
\end{equation}
\end{proposition}
\begin{proof} Most of what is said in the statement is evident, the essential
part being the proof of  \eqref{w bar uno}. This can actually be deduced from a more general
isomorphism between $L_\Pi(A,n+1)$ and the nerve of certain higher dimensional groupoid which is
actually equal to $\W_1\bpar{K(\textcolor{black}{\tilde{A}_1},n)}$ (see
\textcolor{black}{\cite{Garcia2003}}, Prop. 2.4.9). Here we will not make use of these general
facts but will indicate just how to verify the equality at the level of simplices so as to
justify the dimensional jump, and will not bore the reader with the tedious verification of the
equality of faces and degeneracies. Since the $m$-simplices in both simplicial sets are open
horns for every $m>n+2$, its is sufficient to prove
$${\W_1\bpar{K(\textcolor{black}{\tilde{A}_1},n)}}_m = {L_\Pi(A,n+1)}_m,\quad
m=0,\dots,n+2.$$ In the right hand side we have,
$$ {L_\Pi(A,n+1)}_m =
\begin{cases}
\ner (\Pi)_m &\text{if}\quad m\leq n\\
\\
\displaystyle\coprod_{\xi\in\ner (\Pi)_{n+1}} A(x_0^\xi) &\text{if}\quad m = n+1,
\\
\displaystyle\coprod_{\xi\in\ner (\Pi)_{n+2}} A(x_0^\xi)^{n+2} &\text{if}\quad m
= n+2.
\end{cases}
$$ where we used the representation of the generalized Eilenberg MacLane spaces
$L_\Pi(A,n)$ given in \cite{BuFaGa2002}. On the other hand it is clear, using the fact that
$K(\textcolor{black}{\tilde{A}_1},n)_m = (\Pi\xto{1_\Pi}\Pi)\in\Gpd/\Pi, m=0,\dots,n-1$, and the
definition of the $m$-simplices of $\W$ given in \cite{DwKa1984}, that for $m<n$,
$\W_1\bpar{K(\textcolor{black}{\tilde{A}_1},n)}_m = \ner (\Pi)_m.$ By the same token, an
$n$-simplex in $\W_1 \bpar{K(\tilde{A}_1,n)}$ is a sequence $(z_n\xto{u_{n-1}}z_{n-1}\to \cdots
\to z_1\xto{u_0} z_0)$ in $\Pi$ and we get again
$\W_1\bpar{K(\textcolor{black}{\tilde{A}_1},n)}_n = \ner (\Pi)_n.$ Let's now consider an $n+1$
simplex in $\W_1\bpar{K(\textcolor{black}{\tilde{A}_1},n)}$, that is, a sequence $\xi =
(z_{n+1}\xto{(u_n,a)}z_n\xto{u_{n-1}}z_{n-1}\to \cdots \to z_1\xto{u_0} z_0)$ where the $u_i$ are
arrows in $\Pi$ and furthermore $z_{n+1}\xto{(u_n,a)}z_n$ is an arrow in $\semi\Pi A$. This means
that $a$ is an arbitrary element in $A(z_{n}) \cong A(z_{n+1}) = A(x_0^\xi)$ (using the
isomorphism $A(u_n^{-1}):A(z_n) \to A(z_{n+1})$). Thus,
$\W_1\bpar{K(\textcolor{black}{\tilde{A}_1},n)}_{n+1} = \coprod_{\xi\in\ner (\Pi)_{n+1}}
A(x_0^\xi) = {L_\Pi(A,n+1)}_{n+1}$. Finally, let's consider an $n+2$ simplex in
$\W_1\bpar{K(\textcolor{black}{\tilde{A}_1},n)}$, that is, a sequence $\xi =
(z_{n+2}\xto{(u_{n+1},\textcolor{black}{\alpha})}z_{n+1}\xto{(u_{n},a)}z_{n}\to \cdots \to
z_1\xto{u_0} z_0)$ where the $u_i$ are arrows in $\Pi$ and furthermore $z_{n+1}\xto{(u_n,a)}z_n$
and $z_{n+2}\xto{(u_{n+1},\alpha)}z_{n+1}$ are arrows in $\semi\Pi A$ and $\big({\semi\Pi
A}\big)^{n+1}$ respectively. This means that $a\in A(z_n)$ and $\alpha=(a_1,\dots,a_{n+1})\in
{A(z_{n+1})}^{n+1}$, while
$$\xi' = (z_{n+2}\xto{u_{n+1}}z_{n+1}\xto{u_n}z_n\to \cdots \to z_1\xto{u_0}
z_0) \in
\ner (\Pi)_{n+2}.$$ Thus, we get an $n+2$ simplex
$(\xi',\alpha')$ in $L_\Pi(A,n+1)$ with $\alpha' = (a'_0, \dots, a'_{n+1}) \in
{A(z_{n+2})}^{n+2}$ where
\begin{align*} a'_0 &= A(u_{n}u_{n+1})^{-1}(a)\\ a'_1 &= A(u_{n+1})^{-1}(a_1)\\
    &\dots\\ a'_n &= A(u_{n+1})^{-1}(a_n)\\ a'_{n+1} &=
A(u_{n+1})^{-1}\Big(\sum_{i=1}^{n+1}(-1)^{n+1-i} a_{i}\Big).
\end{align*} It is now straightforward to verify that the correspondences between
simplices in $L_\Pi(A,n+1)$ and in $\W_1\bpar{K(\textcolor{black}{\tilde{A}_1},n)}$ we have
established is a $(n+2)$-truncated bijective simplicial map whose $(n+2)$-component satisfies the
cocycle condition, thus determining a simplicial isomorphism.
\end{proof}

A simple extension of Theorem 3.3 in \cite{DwKa1984} yields without difficulty
the following

\begin{proposition}\label{dwka 3.3 extended} The functor $\W_1$ in \eqref{wbar}
preserves fibrations and weak equivalences; it has a left adjoint $\f_1$ which also preserves
fibrations and weak equivalences and for every pair of objects $X\in\scrs{1}/\Pi$,
$Y\in\dsimpl{\sets}/\ner(\Pi)$ a map $Y\to\W_1X$ is a weak equivalence if and only if its adjoint
$\f_1Y\to X$ is a weak equivalence.
\end{proposition}

It follows from this that the adjunction goes through to the corresponding homotopy categories
and as a consequence the set of classes of homotopic maps $Y\to\W_1X$ in
$\dsimpl{\sets}/\ner(\Pi)$ is bijective with the set of classes of homotopic maps $\f_1Y\to X$ in
$\scrs{1}/\Pi$. Taking $X=K(\textcolor{black}{\tilde{A}_1},m)$ and $Y = \ner(\cbc)$, we have

\begin{corollary} \label{cor of dwka 3.3 extended} For any crossed complex
$\cbc$, let $\Pi=\pi_1(\cbc)$ be its fundamental groupoid and let $A:\Pi\to\Ab$ be any
$\Pi$-module. Then, if $\textcolor{black}{\tilde{A}_1} = \semi\Pi A$ is regarded as an abelian
group object in $\Gpd/\Pi$,
$$
\scoH^{m+1}\bpar{\cbc,A} = \big[\f_1\ner(\cbc),
K(\textcolor{black}{\tilde{A}_1},m)\big]_{\scrs{1}/\Pi}.
$$
\end{corollary}

In higher dimensions we will avoid trying to generalize Proposition \ref{dwka 3.3
extended}, proving instead directly the higher dimensional analog of Corollary
\ref{cor of dwka 3.3 extended}.

\subsubsection{The lifting to $\scrs 2$.}

We now define the functor $\W_2:\scrs{2} \to \scrs{1}$.

Let $\sgd \subseteq \ddsimpl{\Gpd}$ be the full subcategory of double simplicial groupoids
determined by those simplicial groupoids all whose vertical and horizontal faces and
degenerations are (functors of groupoids which are) the identity on objects. Equivalently, \sgd\
is the category of groupoids enriched in double simplicial sets. We first notice that the
Artin-Mazur diagonal functor, $\W_{\!\text{A-M}}$, takes object in \sgd\ to objects in
$\scrs{1}$,
$$
\xymatrix{\sgd \ar@{-->}[r]^{\W} \ar@{^{(}->}[d] & \scrs{1}\ar@{^{(}->}[d] \\
\Gpd^{\op{\sdel}\times\op{\sdel}} \ar[r]_{\W_{\!\text{A-M}}} & \Gpd^{\op{\sdel}}
\; .}
$$ On the other hand, there is an isomorphism between the category of crossed
modules and the subcategory $(\scrs{1})_2 \subseteq \scrs{1}$ of simplicial groupoids with
trivial Moore complex in dimensions $\geq 2$ (Corollary 3.1.10 in \cite{Garcia2003}) which,
composed with the inclusion $(\scrs{1})_2\hookrightarrow\scrs{1}$ induces a functor $\dsimpl{\xm}
\rightarrow \dsimpl{\scrs{1}}$ whose restriction to $\scrs{2}$ takes its values in \sgd,
$$
\xymatrix{\scrs{2} \ar@{-->}[r]^\Theta \ar@{^{(}->}[d] & \sgd \ar@{^{(}->}[d] \\
\xm^{\op{\sdel}} \ar[r] & (\scrs{1})^{\op{\sdel}}\; .}
$$ We define $\W_2$ as the composition
$$
\W_2:\scrs{2}\xrightarrow{\Theta} \sgd \xrightarrow{\W} \scrs{1}.
$$ We next describe the action of $\W_2$ on objects. Let
$(\bg,\Sigma)\in\scrs{2}$ so that for each $n\geq0$, $\Sigma_n = (\bg,C_n,\delta_n)$ is a
$\bg$-crossed module. The simplicial groupoid $\W_2(\bg,\Sigma)$ has the same object as \bg. Its
groupoid of 0-simplices is $ \W_2(\bg,\Sigma)_0 = \bg$; its groupoid of 1-simplices can be
identified with $ \W_2(\bg,\Sigma)_1 = \semi \bg C_0$ with face and degeneration morphisms given
by the following formulas:
\begin{equation*}
\begin{split} d_0(f,a_0)& = (\delta_0)_y(a_0)f, \\ d_1(f,a_0)& =d_1^v(f)= f, \\
s_0(f)& = (f,0_{C_0(y)}),
\end{split}
\end{equation*}

In general, for $n\geq 2$, the set of arrows from $x$ to
$y$ in the groupoid of $n$-simplices of $\W_2(\bg,\Sigma)$ is given by
\begin{equation}
\label{formula w bar 2} Hom_{\W_2(\bg,\Sigma)_n}(x,y) = Hom_\bg(x,y) \times
C_0(y)
\times \ldots
\times C_{n-1}(y).
\end{equation} We therefore denote
$$
\W_2(\bg,\Sigma)_n= \bg * C_0 *\cdots*C_{n-1}
$$ the groupoid of $n$-simplices of $\W_2(\bg,\Sigma)$. The composition in this
groupoid is given by the  formula:
\begin{multline*}
 (g,b_0,b_1,\dots,b_{n-1})(f,a_0,a_1,\dots,a_{n-1}) \\ =
(gf,\,b_0+{^{(\delta_1)_z(b_1)\dots (\delta_{n-1})_z(b_{n-1})g}a_0},\dots\\
\dots,b_{n-2}+{^{(\delta_{n-1})_z(b_{n-1})g}a_{n-2}}, b_{n-1}+{^ga_{n-1}}).
\end{multline*}
Note that again this groupoid is a kind of semidirect product.
The face and degeneration operators
$$ (\W_2(\bg,\Sigma))_{n+1} \xleftarrow{s_j} (\W_2(\bg,\Sigma))_n
\xrightarrow{d_i} (\W_2(\bg,\Sigma))_{n-1}
$$ are given by the formulas:
\begin{equation}
\label{cydw2}
\begin{split} d_0(f,a_0,a_1,\dots,a_{n-1}) & =
((\delta_{n-1})_y(a_{n-1})f,a_0,\dots,a_{n-2}), \\
d_i(f,a_0,a_1,\dots,a_{n-1}) & =
(f,a_0,\dots\\
&\phantom{mmn}\dots,a_{n-i-2},a_{n-i-1}+d_0a_{n-i},d_1a_{n-i+1},\dots
,d_{i-1}a_{n-1}),\\
d_n(f,a_0,a_1,\dots,a_{n-1}) & = (f,d_1a_1,\dots,d_{n-1}a_{n-1}),\\
s_j(f,a_0,a_1,\dots,a_{n-1}) & =
(f,a_0,\dots,a_{n-j-1},0_{C_{n-j}(y)},s_0a_{n-j},\dots,s_{j-1}a_{n-1}). \\
\end{split}
\end{equation}

\noindent{\bfseries Example:} If $\Pi$ is a groupoid regarded as a discrete
constant simplicial crossed module then $\W_2(\Pi)$ is equal to $\Pi$ regarded
as a constant simplicial groupoid.

\bigskip

If  $\Pi$ is a groupoid, as a consequence of $\W_2(\Pi) = \Pi$, the functor
$\W_2:\scrs{2} \to
\scrs{1}$ induces a functor $\W_2:\scrs{2}/\Pi \to \scrs{1}/\Pi$. Let $A:\Pi \to
\Ab$ be a
$\Pi$-module. If, as before, for $n\geq1$, $\textcolor{black}{\tilde{A}_n} =
\ins_n(A)$ is regarded as an abelian group object in $\crs_n/\Pi$, we have:

\begin{proposition}
\label{case n eq 2 of the wbars compose to the hocol}
$$
\W_2\bpar{K(\textcolor{black}{\tilde{A}_2}, m)} =
K(\textcolor{black}{\tilde{A}_1},m+1).
$$
\end{proposition}

\begin{proof}
$K(\textcolor{black}{\tilde{A}_2}, m)$ is the simplicial crossed module
$(\Pi,\Sigma)$ where
$\Sigma_n=(\Pi,C_n,0)$ with
$$C_n(y) = \begin{cases} 0 &\text{if }n<m,\\ A(y) &\text{if }n = m,\\ A(y)^{m+1}
&\text{if }n = m+1.
\end{cases}$$ Thus, formula \eqref{formula w bar 2} yields in this case,
$$
\W_2\bpar{K(\textcolor{black}{\tilde{A}_2}, m)}_n(x,y) =
\begin{cases}
\Pi(x,y) &\text{if }n \leq m,\\
\Pi(x,y)\times A(y) &\text{if }n = m+1,\\
\Pi(x,y)\times A(y)^{m+2} &\text{if }n = m+2.
\end{cases}
$$ This is the same one obtains for $K(\textcolor{black}{\tilde{A}_1}, m+1)$.
\textcolor{black}{The details are in \cite{Garcia2003} Cor. 4.3.8.}
\end{proof}

We define now a left adjoint to $\W_2$, which is similar to the ``loop groupoid''
functor $\f_1= G$ defined in \cite{DwKa1984}.

Let $\Sigma_n$ be the groupoid of $n$-simplices of the simplicial groupoid
$\Sigma$. We have a simplicial diagram of pre-crossed modules which in dimension
$n$ has the pre-crossed module
$(\Sigma_0,K_{n-1},\delta)$ which, furthermore has a split augmentation  by
$(\Sigma_0, K_0, \delta)$ where for $n \geq 0$. $K_n$ denotes the $\Sigma_0$-group
associated to the totally disconnected groupoid defined by
$\ker(d_1d_2\cdots d_{n+1}:\Sigma_{n+1}
\rightarrow
\Sigma_0)$ with action given by conjugation via $s_ns_{n-1}\cdots s_0:\Sigma_0
\rightarrow
\Sigma_{n+1}$, and $\delta$ is the natural transformation whose $x$-component,
for $x \in
\obj(\Sigma_0)$, is $$\delta_x(u) = \begin{cases}d_0(u), \text{ if } u \in
K_0(x)\\ d_0d_2\cdots d_{n+1}(u) \text{ if }u \in K_n(x)\ \text{with}\ n\geq
1.\end{cases}$$  Note that the  face and degeneration operators
$$ (\Sigma_0, K_{n+1}, \delta)\xleftarrow{\sigma_j} (\Sigma_0, K_n,
\delta)\xrightarrow{\delta_i} (\Sigma_0, K_{n-1}, \delta)
$$ for $1 \leq i \leq n$ and $0 \leq j \leq n$, and also the augmentation
$\sigma_0:(\Sigma_0, K_n, \delta)
\rightarrow (\Sigma_0, K_{n+1}, \delta)$ for all $n \geq 0$, are given by
restrictions of the faces
$d_{i+1}:\Sigma_{n+1}
\rightarrow \Sigma_n$ and of the degeneracies $s_{j+1}:\Sigma_{n+1}
\rightarrow \Sigma_{n+2}$ of the simplicial groupoid $\Sigma$.

From this we build an augmented split simplicial crossed module by
factoring out in each $\Sigma_0$-module, the Peiffer elements as
well as those who are images by $s_0$. Such quotient determines,
in each dimension, a crossed module
$(\Sigma_0,\tilde{K}_n,\delta)$ and the face and degeneration
operators go well with the quotients. In order to obtain the
simplicial crossed module $\f_2(\Sigma)$ we just need to add in
each dimension a new  morphism of crossed modules which will be
given by the morphism of $\Sigma_0$-groups $\delta_0:\tilde{K}_n
\rightarrow \tilde{K}_{n-1}$, in turn determined by the natural
transformation $[d_0,d_1]:K_n \rightarrow K_{n-1}$ whose component
on an object $x \in \obj(\Sigma_0)$ is given by
$$ [d_0,d_1]_x(u) = (d_1)_x(u) (d_0)_x(u)^{-1}(s_0d_1d_0)_x(u)
$$ for each $u \in K_n(x)$. The functor so defined is left adjoint to $\W_2$ (see
\cite{Garcia2003} Prop. 4.3.9).

If we regard now a groupoid $\Pi$ first as a constant simplicial groupoid and
then as a constant simplicial crossed module, it is easy to check that
$\f_2(\Pi)=\Pi$. In fact, we still have an adjoint pair of functors
$$\f_2\dashv\wbar[2], \quad
\xymatrix@C=2.2pc@R=1.2pc{ {\vphantom{\Gr^\bg}}\scrs{2}/\Pi
\ar@<-.6ex>[r]_-{\wbar[2]}
&{\vphantom{\pxm[\bg]}}\scrs{1}/\Pi\ar@<-.4ex>[l]_-{\f_2}
  } .$$ We next show that in certain cases the functors $\f_2$ and
$\W_2$ preserve homotopy classes. First we note that in order that two morphisms
$(F,\boldsymbol{\alpha}),(G,\boldsymbol{\beta}):(\bg, \Sigma) \rightarrow
(\bg',\Sigma')$ in
$\scrs{2}$ be homotopic, both must have the same functor at the level of base
groupoids, that is,
$F=G$ and also that if
$\hbar=(F,\bh)$ is a homotopy between them with
$$
\bh=\{h_j^n;\; 0\leq j\leq
n\}:\boldsymbol{\alpha}\rightsquigarrow\boldsymbol{\beta}:
\Sigma\rightarrow\Sigma'
$$ where $h_j^n=(F,\gamma_j^n)$, the homotopy identities for $\hbar$ are
essentially the homotopy identities for the homotopy
$\boldsymbol{\gamma}=\{\gamma_j^n:C_n\rightarrow C'_{n+1}F\}$ between
simplicial morphisms of
$\bg$-groups.

\begin{lemma} The functor $\W_2$ preserves homotopy classes of simplicial
morphisms.
\end{lemma}
\begin{proof} Given a homotopy $\hbar=(F,\bh)$ as before, if we denote
$\mathbf{F}=\W_2(F,\boldsymbol{\alpha})$ and $\mathbf{F}'=\W_2(F,\boldsymbol{\beta})$, then the
components of $\mathbf{F}$ and $\mathbf{F}'$ in dimension $0$ are given by the functor
$F_0=F'_0=F$ and in dimension $n$ by the functors
$$ F_n, F'_n: \bg*C_0*\cdots *C_{n-1} \rightarrow \bg'*C'_0*\cdots *C'_{n-1}.
$$ These functors act on objects as the functor $F$ and, on an arrow
$(f,a_0,\dots,a_{n-1})$ with codomain $y$, they act thus:
$$
\begin{array}{c} F_n(f,a_0,\dots,a_{n-1}) =
\big(F(f),(\alpha_0)_y(a_0),\dots,(\alpha_{n-1})_y(a_{n-1})\big),\\ [1pc]
F'_n(f,a_0,\dots,a_{n-1}) = \big(F(f),(\beta_0)_y(a_0),\dots,(\beta_{n-1})_y(a_{n-1})\big).
\end{array}
$$ We have a homotopy $\mathbf{H}:\mathbf{F}\rightsquigarrow \mathbf{F}':
\W_2(\bg,\Sigma)
\rightarrow \W_2(\bg',\Sigma')$ with
$$
\mathbf{H}=\{H_j^n:\bg*C_0*\cdots *C_{n-1} \rightarrow \bg'*C'_0*\cdots *C'_n;\; 0
\leq j \leq n
\}.
$$ where  $H_j^n:\bg*C_0*\cdots *C_{n-1} \rightarrow \bg'*C'_0*\cdots *C'_n$, for
$0\leq j \leq n$, acts on objects as the functor $F$ and, on each arrow
$(f,a_0,\dots,a_{n-1})$ with codomain $y$, it acts thus:
\begin{multline*}
H_j^n(f,a_0,\dots,a_{n-1})=
\big(F(f),(\alpha_0)_y(a_0),\dots\\
\dots,(\alpha_{n-j-1})_y(a_{n-j-1}),0,
(\gamma_0^{n-j})_y(a_{n-j}),\dots, (\gamma_{j-1}^{n-1})_y(a_{n-1})\big).
\end{multline*}
It is easy to check that the homotopy identities for $\mathbf{H}$
follow from the corresponding identities for $\boldsymbol{\gamma}$. The details
are in
\cite{Garcia2003}.
\end{proof}

The functor $\f_2$ does not behave the same way as $\W_2$ with respect to
homotopies. However, in certain cases $\f_2$ does take homotopic morphisms to
homotopic morphisms. One of these cases is the following:
\begin{lemma} Let $\mathbf{F},\mathbf{F}':\Sigma \rightarrow \Sigma'$ be two
morphisms in
$\scrs{1}$ such that
$F_0 = F'_0$ and let $\mathbf{H}$ be a homotopy between them, such that
$$ H_0^0 = s_0F_0 \qquad\mbox{and}\qquad H_i^j = s_iF_j \mbox{\; if \;} i<j
$$ (note that $H_j^j$ is arbitrary for $j>0$), then the morphisms of simplicial
crossed modules $\f_2(\mathbf{F})$ and $\f_2(\mathbf{F}')$ in $\scrs{2}$ are
homotopic.
\end{lemma}

The two previous lemmas can easily be proved for slice categories. We have:

\begin{lemma}
\label{w2 conserva homo} For each groupoid $\Pi$, the functor $\W_2:\scrs{2}/\Pi
\rightarrow
\scrs{1}/\Pi$ preserves homotopy classes of simplicial morphisms in the
corresponding slice categories.
\end{lemma}
\begin{lemma}
\label{f2 conserva homo} Let $\Pi$ be a groupoid, $\mathbf{F}$ and $\mathbf{F}'$
two morphisms in
$\scrs{1}/\Pi$
$$
\xymatrix{ \Sigma \ar@<0.6ex>[rr]^{\mathbf{F}} \ar@<-0.6ex>[rr]_{\mathbf{F}'}
\ar[dr] & & \Sigma'
\ar[dl] \\ &\Pi& }
$$ such that $F_0 = F'_0$ and let $\mathbf{H}$ be a homotopy between them in
$\scrs{1}/\Pi$
$$ H_0^0 = s_0F_0 \qquad\mbox{and}\qquad H_i^j = s_iF_j, \mbox{ if }\ i<j.
$$ Then the morphisms of simplicial crossed modules $\f_2(\mathbf{F})$ and
$\f_2(\mathbf{F}')$ in the slice category $\scrs{2}/\Pi$ are homotopic.
\end{lemma}

Note that if in Lemma \ref{f2 conserva homo} we take $\Sigma'=K(\tilde{A}_1,n) \in \scrs{1}/\Pi$
for any $\Pi$-module $A$, then for any two homotopic morphisms in $\scrs{1}/\Pi$ with codomain
$K(\tilde{A}_1,n)$ there is a homotopy $\mathbf{H}$ in the hypothesis of the said lemma and
therefore we can conclude that the functor $\f_2$ preserves homotopy classes of morphisms with
codomain $K(\tilde{A}_1,n)$.

From the preceding reasoning it follows immediately,

\begin{proposition} Let $\Pi$ be a groupoid, $X$  a simplicial groupoid above
$\Pi$, and let $A$ be $\Pi$-module. Then the adjunction $\f_2\dashv\wbar[2]$,
$$
\xymatrix@C=2.2pc@R=1.2pc{ {\vphantom{\Gr^\bg}}\scrs{2}/\Pi
\ar@<-.6ex>[r]_-{\wbar[2]}
&{\vphantom{\pxm[\bg]}}\scrs{1}/\Pi\ar@<-.4ex>[l]_-{\f_2}
  }
$$ induces an isomorphism in homotopy classes:
\begin{equation}
\label{eq iso in homotopy}
\left[\f_2(X), K(\abgp A2,n) \right]_{\scrs{2}/\Pi} \cong
\left[X , \wbar[2]\bpar{K(\abgp A2,n)} \right]_{\scrs{1}/\Pi} .
\end{equation}
\end{proposition}

\bigskip

\subsubsection{The lifting to $\scrs n$ for $n \geq 3$.}

For $n\geq3$ each object of $\scrs{n}$ can be represented by a diagram of this
form:
\begin{equation}\label{scrsn obj}
\xy
   (-10,0)*+{\ldots}="z",
   (0,0)*+{\c_n^i}="a",
   (20,0)*+{\c_n^{i-1}}="b",
   (30,0)*+{\ldots}="c",
   (40,0)*+{\c_n^1}="d",
   (60,0)*+{\c_n^0\,.}="e",
   (30,-10)*+{\c_{n-1}}="g",
   (30,-20)*+{\c_{n-2}}="h",
   (30,-30)*+{\vdots}="i",
   (30,-40)*+{\c_2}="j",
   (30,-50)*+{\termxm[\bg]}="k",
\ar @/^-4ex/ "b";"a" |{s_0}
\ar @/^-5ex/@{..} "b";"a"
\ar @/^-6ex/ "b";"a" |{s_{i-1}}
\ar @{->}^{d_i} "a";"b" <3pt>
\ar @{..} "a";"b"
\ar @{->}_{d_0} "a";"b" <-3pt>
\ar @/^-4ex/ "e";"d" |{s_0}
\ar @{->}^{d_1} "d";"e" <3pt>
\ar @{->}_{d_0} "d";"e" <-3pt>
\ar @/^-3ex/ "a";"g"_{\partial_n^i}
\ar @{->}_{\partial_n^{i-1}} "b";"g"
\ar @{->}^{\partial_n^1} "d";"g"
\ar @/^3ex/ "e";"g"^{\partial_n^0}
\ar @{->}^{\partial_{n-1}} "g";"h"
\ar @{->} "h";"i"
\ar @{->} "i";"j"
\ar @{->}^{\partial_2} "j";"k"
\endxy
\end{equation} Giving such an object is equivalent to giving it's ``head'':
\begin{equation}\label{scrsn head}
\xy
   (-20,0)*+{\bc_n:}="w",
   (-10,0)*+{\ldots}="z",
   (0,0)*+{\c_n^i}="a",
   (20,0)*+{\c_n^{i-1}}="b",
   (30,0)*+{\ldots}="c",
   (40,0)*+{\c_n^1}="d",
   (60,0)*+{\c_n^0}="e",
\ar @/^-4ex/ "b";"a" |{s_0}
\ar @/^-5ex/@{..} "b";"a"
\ar @/^-6ex/ "b";"a" |{s_{i-1}}
\ar @{->}^{d_i} "a";"b" <3pt>
\ar @{..} "a";"b"
\ar @{->}_{d_0} "a";"b" <-3pt>
\ar @/^-4ex/ "e";"d" |{s_0}
\ar @{->}^{d_1} "d";"e" <3pt>
\ar @{->}_{d_0} "d";"e" <-3pt>
\endxy
\end{equation} (a simplicial complex of $\bg$-modules, where $\bg$ is the base
groupoid of all the involved crossed complexes), its ``tail''
$$
\cbc=(\c_{n-1} \xto{\partial_{n-1}} \c_{n-2} \rightarrow \ldots \to \c_2
\xto{\partial_2}
\termxm[\bg])
$$ (an $(n-1)$-crossed complex), and an augmentation
$\bc_n\xrightarrow{\partial_n^0} C_{n-1}$ of
$\bc_n$ over $C_{n-1}=\techo_{n-1}(\cbc)$ such that the compositions
$$ C_n^0\xrightarrow{\partial_n^0}C_{n-1}\xrightarrow{\partial_{n-1}}
C_{n-2}\qquad\mbox{and}\qquad
\hat{C}_2\xrightarrow{\partial_2}\g\xrightarrow{C_n^0}\ab
$$ are trivial.

We will represent the above  $n$-crossed complex by the pair $(\cbc,\bc_n)$ or
the triple
$(\bg,\cbc,\bc_n)$ in case we want to make explicit the base groupoid.

A map from $(\cbc, \bc_n)$ to $(\cbc', \bc_n')$ in $\scrs n$ is a pair
$(\bf,\boldsymbol{\alpha})$ where
$\bf:\cbc \rightarrow \cbc'$ is a morphism of $(n-1)$-crossed complexes with $F=
\base (\bf)$ and $\boldsymbol{\alpha}:\bc_n \rightarrow F^\ast\,\bc_n'$
is a  simplicial map of
$\bg$-modules where $\bg=\base(\cbc)$, $\bg'=\base(\cbc')$, and
$F^\ast:\ab^{\bg'} \to \ab^\bg$ is the functor induced by $F=\base(\bf)$.

It is possible to give a definition of $\W_n$ in terms of the Artin-Mazur
diagonal as for the case $n=2$. But for $n>2$ it is also possible to give a
direct description without going through double simplicial $n$-crossed complexes:
$$
\W_n(\cbc, \bc_n) = (T_{n-2}(\cbc), \bc_{n-1})
$$ where $\bc_{n-1}$ is the simplicial $\bg$-module \big(where $\bg=\base(\cbc)$\big)
augmented over
$C_{n-2} = \techo_{n-2}(\cbc)$ given by
$$ C^0_{n-1}=C_{n-1}\qquad\mbox{and}\qquad C^i_{n-1}= C_n^{i-1} \oplus \cdots
\oplus C_n^0 \oplus C_{n-1},\qquad\mbox{if}\; 1\leq i,
$$ where $\c_n^i =(\bg,C_n^i,0)$ for $i \geq 0$ and with augmentation
$\partial_{n-1}^0=\partial_{n-1}:C_{n-1}\rightarrow C_{n-2}$ and with face and
degeneration operators
\begin{align*}
\label{cydwn} (d_0)_x(u_{i-1},\dots,u_0,u) & =(u_{i-2},\dots,u_0,
\partial_n^{i-1}(u_{i-1})+u)\\ (d_i)_x(u_{i-1},\dots,u_0,u) & =
(d_{i-1}u_{i-1},\dots,d_1u_1, u)\\ (d_j)_x(u_{i-1},\dots,u_0,u) & =
(d_{i-1}u_{i-1},\dots,d_1u_{i-j+1},d_0u_{i-j}+u_{i-j-1},\dots,u_0,u)\tag{if  $1\leq
j < i$}\\ (s_j)_x(u_{i-1},\dots,u_0,u) & =
(s_{j-1}u_{i-1},\dots,s_0u_{i-j},0,u_{i-j-1},\dots,u_0, u)\tag{if   $0
\leq j \leq i$}.
\end{align*} It is obvious how $\W_n$ acts on morphisms.

If  a groupoid  $\bg$ is regarded as the simplicial $n$-crossed complex
$(\bg,0)$, we have
$$
\W_n(\bg)=\W_n(\bg,0)=(\bg,0)=\bg
$$ As a consequence, $\W_n$ induces a functor between slice categories,
$$
\W_n: \scrs n/\bg \to \scrs{n-1}/\bg
$$
\begin{lemma} If $(\cbc,\bc_n)\in\scrs n$ is such that $\bc_n$ has trivial Moore
complex in dimensions  $\geq m$,  then upon applying $\W_n$, we get $\W_n(\cbc,
\bc_n) = (T_{n-2}(\cbc), \bc_{n-1})$, where
$\bc_{n-1}$  again has trivial Moore complex in dimensions $\geq m+1$.
\end{lemma}

\noindent{\bfseries Example:} Let $\Pi$ be a groupoid and $A:\Pi
\rightarrow
\ab$ a
$\Pi$-module, for each $m \geq 1$ the triple $(\Pi,\Pi,K(A, m))$ determines a
$n$-simplicial crossed complex in $\scrs n$. By the previous lemma,
$$
\W_n\big(\Pi, \Pi, K(A, m)\big) \cong (\Pi,\Pi,K(A,m+1)).
$$ On the other hand,
$$ K(\tilde\Pi_n,m) = \big(\Pi, \Pi, K(A, m)\big)
$$ and therefore we have,
\begin{proposition}
\label{the wbars compose to the hocol} For every $n \geq 1$ ,
$$
\W_n\big(K(\tilde\Pi_n,m)\big) = K(\tilde\Pi_{n-1},m+1).
$$
\end{proposition}

We define now the left adjoint to $\W_n$ on objects. Given $(\bg,\cbc,\bc_{n-1})
\in \scrs{n-1}$, the $\bg$-modules
$$
\begin{array}{lc} K_0=\Ker(d_1:C_{n-1}^1 \rightarrow C_{n-1}^0), &  \\ [1pc] K_i
=
\Ker(d_1d_2\cdots d_{i+1}:C_{n-1}^{i+1} \rightarrow C_{n-1}^0), & i \geq 1.
\end{array}
$$ together with the  operators $d_j:K_i \rightarrow K_{i-1}$ and
$s_{j-1}:K_{i-1} \rightarrow K_i$ induced by the face operators $d_j:C_{n-1}^{i+1} \rightarrow
C_{n-1}^i$ and degeneracies $s_{j-1}:C_{n-1}^i \rightarrow C_{n-1}^{i+1}$, for $2 \leq j \leq
i+1$ define an augmented split simplicial complex of $\bg$-modules. Factoring out by the
$s_0:C_{n-1}^i \rightarrow C_{n-1}^{i+1}$-image of $K_{i-1}$ we get $\bg$-modules
$$
\widetilde{K}_i= K_i/ s_0(K_{i-1})
$$ which again determine an  augmented split simplicial complex of $\bg$-modules
with face operators
$\delta_j$ and degeneracies $\sigma_j$ induced by the corresponding quotients by
the operators $d_{j+1}$ and
$s_{j+1}$. Finally, the natural transformation whose component in an object
$x\in \obj(\bg)$ is given by:
$$ (\delta_0)_x(u) = (d_1)_x(u)\, (d_0)_x(u)^{-1}(s_0d_1d_0)_x(u),
$$ for each $u \in K_i(x)$, determines another, also be denoted $\delta_0:
\widetilde{K}_i
\rightarrow \widetilde{K}_{i-1}$, for $i>0$, which together with the previous
diagram provides us with a simplicial complex  of
$\bg$-modules which will be denoted $\widetilde{\mathbf{K}}_{n}$. Furthermore,
the restriction to $K_0$ of the face operator
$d_0:C_{n-1}^1
\rightarrow C_{n-1}^0$  induces a morphism
$\partial_n^0:\widetilde{K}_0 \rightarrow C_{n-1}^0$ such that
$\widetilde{\mathbf{K}}_n
\xrightarrow{\partial_n^0} C_{n-1}^0$ is an augmented simplicial $\bg$-module.
Thus, we define
$$
\f_n(\bg,\cbc,\bc_{n-1}) = (\bg,\cbc_{n-1}^0,\widetilde{\mathbf{K}}_n),
$$ where $\cbc_{n-1}^0$ is the $(n-1)$-crossed complex whose $n-1$ truncation is
$\cbc$
($T_{n-1}(\cbc_{n-1}^0) =\cbc$), and is  such that
$\techo_{n-1}(\cbc_{n-1}^0) = C_{n-1}^0$.

It is now a routine exercise to define this functor on arrows and to verify
that it is left adjoint to $\W_n$.

\begin{note}
For every $n > 3$, the functor $\W_n:\scrs{n} \to \scrs{n-1}$ is an equivalence
of categories whose inverse is $\f_n$, and as a consequence, for $n>3$ the categories
$\scrs{n}$ are equivalent to $\scrs{3}$.
\end{note}

Let $\Pi$ be a groupoid, regarded as a \ncc[1]. Evidently $\Pi$ can be regarded
as a \ncc[k] for any $k$, and we will do so as needed. Furthermore, we can
consider the constant simplicial \cc\
$\cte{\Pi}$ as an object in $\scrs{k}$ for any $k$. Such objects verify
$\W(\cte{\Pi}) =
\cte{\Pi}$ so that the functors $\W_k:\scrs{k} \to \scrs{k-1}$ induce functors
$\W_k:\scrs{k}/\cte{\Pi} \to \scrs{k-1}/\cte{\Pi}$. In fact, for every $k\geq1$,
we still have an adjoint pair of functors
$$
\f_k\dashv\wbar[k], \quad
\xymatrix@C=2.2pc@R=1.2pc{ {\vphantom{\Gr^\bg}}\scrs{k}/\Pi
\ar@<-.6ex>[r]_-{\wbar[k]}
&{\vphantom{\pxm[\bg]}}\scrs{k-1}/\Pi\ar@<-.4ex>[l]_-{\f_k}
  }
$$ We next show that in certain cases the functors $\f_n$ and
$\W_n$ preserve homotopy classes. First we note that  the restriction to $\scrs
n$ of the functor
$(T_{n-1})_\ast:\dsimpl{{\crs_n}} \rightarrow \dsimpl{{\crs_{n-1}}}$ takes
homotopies to homotopies. Thus, a homotopy
$\hbar:(F,\bf,\boldsymbol{\alpha})\rightsquigarrow (F,\bf,\boldsymbol{\beta})$
is given as a pair
$$
\hbar=(\bf,\bh),
$$ with
$$
\bh =\{h_j^i:C_n^i \rightarrow C_n^{\prime\, i+1}F;\; 0\leq j\leq
i\}:\boldsymbol{\alpha}\rightsquigarrow \boldsymbol{\beta}: \bc_n\rightarrow
F^\ast \bc_n'
$$ a homotopy in the category $\dsimpl{(\ab^\bg)}$. Furthermore, the  homotopy
identities for $\hbar$ follow from the  homotopy identities for $\bh$.

We can now prove

\begin{lemma}
\label{wnpch} The functor $\W_n$ preserves homotopy classes of  simplicial
morphisms.
\end{lemma}
\begin{proof} Let $\hbar=(\bf,\bh)$ be a  homotopy as before, and let us put
$\W_n(F,\bf,\boldsymbol{\alpha})= \big(F,T_{n-2}(\bf),\boldsymbol{\alpha}'\big)$
and $\W_n(F, \bf, \boldsymbol{\beta}) = \big(F, T_{n-2}(\bf),
\boldsymbol{\beta}'\big)$, where $\boldsymbol{\alpha}'$ and $\boldsymbol{\beta}'$
are simplicial morphisms of $\bg$-modules given in dimension  $i$ by
$$
\begin{array}{c}
\alpha_{n-1}^{\prime\, i}(u_{i-1},\dots,u_0,u) = (\alpha_n^{i-1}(u_{i-1}),\dots,
\alpha_n^0(u_0),\alpha_{n-1}(u)),\\ [1pc] \beta_{n-1}^{\prime\,
i}(u_{i-1},\dots,u_0,u) =
(\beta_n^{i-1}(u_{i-1}),\dots,\beta_n^0(u_0),\beta_{n-1}(u)),
\end{array}
$$ for each $(u_{i-1},\dots,u_0,u) \in C_n^{i-1}(x)\oplus \cdots  \oplus C_n^0(x)
\oplus C_{n-1}(x)$. Then the homotopy $\hbar'=(T_{n-2}(\bf),\bar{\bh})$ with
$$
\bar{\bh}=\{\bar{h}_j^i;\; 0 \leq j \leq i\}:\boldsymbol{\alpha}' \rightarrow
\boldsymbol{\beta}'
$$ where $\bar{h}_j^i$ is the natural transformation whose $x$-component for $x
\in \obj(\bg)$ acts thus:
\begin{multline*}
\hbar_j^i(u_{i-1},\dots,u_0,u) =
\big(h_{j-1}^{i-1}(u_{i-1}),\dots,h_0^{i-j}(u_{i-j}),0,\alpha_n^{i-j-1}(u_{i-j-1}),\dots
\\
\dots, \alpha_n^0(u_0), \alpha_{n-1}(u)\big),
\end{multline*}
for each $(u_{i-1},\dots,u_0,u) \in C_n^{i-1}(x)\oplus \cdots
\oplus C_n^0(x) \oplus C_{n-1}(x)$. It is immediate to check that the homotopy
identities for
$\bar{\bh}$ follow from the corresponding identities satisfied by $\bh$ (see
\cite{Garcia2003} Lema 4.3.23 for the details).
\end{proof}

The functor $\f_n$ does not behave the same way as $\W_n$ with respect to
homotopies. However, in certain cases $\f_n$ does take homotopic morphisms to
homotopic morphisms. One of these cases is the following:

\begin{lemma} Let $(F,\bf,\boldsymbol{\alpha}),
(F,\bf,\boldsymbol{\beta}):(\bg,\cbc,\bc_{n-1}) \rightarrow (\bg',\cbc',\bc'_{n-1})$ be morphisms
in the category $\scrs{n-1}$ such that $\alpha_{n-1}^0 = \beta_{n-1}^0:C_{n-1}^0 \rightarrow
C_{n-1}^{' 0}F$ and let $\hbar=(\bf,\bh)$ be a homotopy between the above morphisms such that
$$ h_0^0 = s_0\alpha_{n-1}^0 \qquad\mbox{and}\qquad h_i^j = s_i \alpha_{n-1}^j
\mbox{\; if \;} i<j
$$ (note that $h_j^j$ is arbitrary for $j>0$). Then the morphisms of simplicial
$n$-crossed complexes
$\f_n(F,\bf,\boldsymbol{\alpha})$ and $\f_n(F,\bf,\boldsymbol{\beta})$ in
$\scrs{n}$ are homotopic.
\end{lemma}

The two previous lemmas can easily be proved for slice categories. We have:

\begin{lemma}
\label{wn conserva homo} For each groupoid $\Pi$, the functor $\W_n:\scrs{n}/\Pi
\rightarrow
\scrs{n-1}/\Pi$ preserves homotopy classes of simplicial morphisms in the
corresponding slice categories.
\end{lemma}
\begin{lemma}
\label{fn conserva homo} Let $\Pi$ be a groupoid, $(F,\bf,\boldsymbol{\alpha})$
and
$(F,\bf,\boldsymbol{\beta})$ two morphisms in the slice category
$\scrs{n-1}/\Pi$,
$$
\xymatrix{ (\bg, \cbc,\bc_{n-1}) \ar@<0.6ex>[rr]^{(F,\bf, \boldsymbol{\alpha})}
\ar@<-0.6ex>[rr]_{(F,\bf,\boldsymbol{\beta})} \ar[dr] & & (\bg',\cbc',\bc'_{n-1})
\ar[dl] \\ &\Pi & },
$$ such that $\alpha_0 = \beta_0:C_{n-1}^0 \rightarrow C_{n-1}^{\,\prime\, 0}F$,
and let
$\hbar=(\bf,\bh):(F,\bf, \boldsymbol{\alpha}) \rightsquigarrow
(F,\bf,\boldsymbol{\beta})$, with
$$
\bh=\{h_j^i:C_{n-1}^i \rightarrow C_{n-1}^{\, \prime\, i+1}F,\; 0 \leq j \leq
i\},
$$ be a homotopy in
$\scrs{n-1}/\Pi$ from $(F, \bf,\boldsymbol{\alpha})$ to $(F,
\bf,\boldsymbol{\beta})$ such that
$$ h_0^0 = s_0 \alpha_{n-1}^0 \qquad \mbox{and} \qquad h_i^j= s_i \alpha_{n-1}^j
\; \mbox{if $i < j$}.
$$ Then the morphisms $\f_n(F,\bf, \boldsymbol{\alpha})$ and $\f_n(F,\bf,
\boldsymbol{\beta})$ in the slice category $\scrs{n}/\Pi$ are homotopic.
\end{lemma}

Note that if $A$ is  any $\Pi$-module and in Lemma \ref{fn conserva homo} we take
$$(\bg',\cbc',\bc'_{n-1}) = K(\tilde{A}_{n-1},m)=(\Pi,\Pi,K(A,m)) \in
\scrs{n-1}/\Pi$$ for $m>0$, then for any two homotopic
morphisms in $\scrs{n-1}/\Pi$ with codomain $K(\tilde{A}_{n-1},m)$ there is a
homotopy
$\hbar$ satisfying the hypothesis of the said lemma and therefore we can
conclude that the functor $\f_n$ preserves homotopy classes of morphisms with
codomain $K(\tilde{A}_{n-1},m)$ and we have,

\begin{proposition} Let $\bg$ be a groupoid, $A$ a $\Pi$-module and
$(\bg,\cbc,\bc_{n-1})$ a simplicial object in
$\scrs{n-1}/\Pi$. Then the adjunction $\f_n\dashv\wbar[n]$,
$$
\xymatrix@C=2.2pc@R=1.2pc{ {\vphantom{\Gr^\bg}}\scrs{k}/\Pi
\ar@<-.6ex>[r]_-{\wbar[n]}
&{\vphantom{\pxm[\bg]}}\scrs{k-1}/\Pi\ar@<-.4ex>[l]_-{\f_n}
  }
$$ induces an isomorphism in homotopy:
\begin{equation}
\label{eq iso in homotopyn} \left[\f_n(\bg,\cbc,\bc_{n-1}), K(\tilde{A}_n,m)
\right]_{\scrs{n}/\Pi} \cong \left[(\bg,\cbc,\bc_{n-1}) ,
\wbar[k](K(\tilde{A}_n,m))
\right]_{\scrs{n-1}/\Pi} .
\end{equation}
\end{proposition}
\begin{proof} The proof follows trivially after the previous observations,
together with Lemma
\ref{wn conserva homo} and the fact that the adjunction isomorphisms for
$\W_n\dashv
\f_n$ are obtained by applying the functors $\W_n$ and $\f_n$, and composing
with the unit and counit of the said adjunction.
\end{proof}

\subsubsection{The representation of singular cohomology of $n$-crossed
complexes as homotopy classes of maps of simplicial
$n$-crossed complexes.\\ }

Let now $${\frak{F}}_n = \f_n\f_{n-1}\cdots \f_2
\f_1.$$ Then, by repeated application of \eqref{eq iso in homotopy} one gets:
$$
\big[{\frak{F}}_n\ner(\cbc), K(\abgp An,m)\big]_{\scrs{n}/\Pi} \xto{\ \
\cong\ \ }
\big[\ner(\cbc), \W_1\cdots\W_n\bpar{K(\abgp An, m)}\big]_{\sset/\ner(\Pi)}.
$$ Combining this isomorphism with Proposition \ref{the wbars compose to the
hocol}, and the generalized Eilenberg-MacLane representation
\eqref{gener eilmacl} of the singular cohomology, we get:

\begin{corollary}
\label{cor nat iso} Let $1\leq m \leq n$. There is a natural isomorphism
$$
\scoH^{n+m}(\cbc,A) \xto{\ \ \cong\ \ } \big[{\frak{F}}_n\ner(\cbc), K(\abgp
An,m)\big]_{\scrs{n}/\Pi}.
$$
\end{corollary}

\subsection{Obtaining the topological invariants from the algebraic ones}
\label{sec the top invars}

The main tool we use in this section is Duskin's representation theorem,
which in our context can be particularized like this:

\begin{theorem}[Duskin's Representation Theorem]
\label{duskins repr theo}\ \nopagebreak

If $\tilde{A}_n$ is any internal abelian group object in $\scrs{n}/\Pi$,
$$ H^m_{\cotrip[n]}(\cbc, \tilde{A}_n) \cong \big[\cotrsr(\cbc), K(\tilde{A}_n,
m)\big]_{\scrs{n}/\Pi} .
$$
\end{theorem} We use this theorem together with Corollary \ref{cor nat iso} to
prove the following:

\begin{theorem} Let $\cbc\in\crs_n$, $\Pi = \pi_1(\cbc)$, and let $A$ be a
$\Pi$-module, $A:\Pi \to \Ab$, and $\tilde{A}_n$ the abelian group object in
$\scrs{n}/\Pi$ defined in \eqref{abelian group in xm}. For every $m\geq0$ there
is a natural map
$$ H^m_{\cotrip[n]}(\cbc, \abgp An) \xto{\ \ \ } \scoH^{n+m}(\cbc, A).
$$
\end{theorem}

\begin{proof} By Corollary \ref{cor nat iso} and Duskin's representation theorem
(Theorem \ref{duskins repr theo}) it is sufficient to define a map
$$
\alpha : \big[\cotrsr(\cbc), K(\abgc An, m) \big]_{\scrs{n}/\Pi} \xto{\ \ \ }
\big[{\frak{F}}_n\ner(\cbc), K(\abgp An,m)\big]_{\scrs{n}/\Pi}.
$$ For this, in turn, it is sufficient to give a morphism
\begin{equation}
\label{the map} \boldsymbol{\eta} : {\frak{F}}_n\ner(\cbc) \to \cotrsr(\cbc)
\end{equation} in $\scrs{n}$ such that it defines a map in $\scrs{n}/\Pi$. Note
that in that case the induced map of sets
$$
\boldsymbol{\eta}_* : \scrs{n}\bpar{\cotrsr(\cbc), K(\abgc An, m)} \xto{\ \ \ }
\scrs{n}\bpar{{\frak{F}}_n\ner(\cbc), K(\abgp An,m)}
$$ automatically preserves homotopy classes.

To define $\boldsymbol{\eta}$ we first observe that the simplicial object
${\frak{F}}_n\ner(\cbc)$ is free with respect to the cotriple \cotrip[n]
\cite{Garcia2003} and therefore to give $\boldsymbol{\eta}$ it is sufficient to
give a morphism $\eta_{-1}:\pi_0{\frak{F}}_n\ner(\cbc)\rightarrow \cbc$  where
$\pi_0:\scrs n \rightarrow \crs_n$ is the connected components functor, defined
as the left adjoint to the diagonal $\Delta: \crs_n \rightarrow \scrs n$. Then
we consider the two adjunctions
$$
\pi_0\dashv \Delta, \;\everyentry={\vphantom{\Bigg(}} \xymatrix@C=1.5pc@R=1.2pc{
\crs_n \ar@<-.6ex>[r]_-{\Delta} & \scrs{n}\ar@<-.4ex>[l]_-{\pi_0} } \mbox{\ \
and\ \ } {\frak{F}_n}\dashv {\frak W}_n=\W_1\cdots\W_n,\;
\xymatrix@C=1.5pc@R=1.2pc{ \scrs{n} \ar@<-.6ex>[r]_-{{\frak W}_n} &
\scrs{0}\ar@<-.4ex>[l]_-{{\frak F}_n}, }
$$ and observe that the composition $\frak{W}_n \Delta$ is just the nerve
functor $\ner$. Therefore $\pi_0{\frak{F}}_n\ner(\cbc)=\pi_0\frak{F}_n\frak{W}_n
\Delta (\cbc)$ and we can take $\eta_{-1}$ as the $\cbc$-component of the counit
of the adjunction $\pi_0\frak{F}_n \dashv \frak{W}_n \Delta$.
\end{proof}

By taking the case $m=2$ we obtain the desired map
$$ H^2_{\cotrip[n]}\bpar{\cbc, \abgc An} \xto{\ \ \ } \scoH^{n+2}\bpar{\cbc, A}.
$$ This implies that for an arbitrary crossed complex $\cbc$, if we denote
$\tilde\pi^{(n)}_{n+1}(\cbc)$ the abelian group object $\abgc An$ corresponding to the
local coefficient system determined by $A=\tilde\pi_{n+1}(\cbc)$, then we have a
morphism
$$ H^2_{\cotrip[n]}\bpar{P_n(\cbc), \tilde\pi^{(n)}_{n+1}(\cbc)} \xto{\ \ \ }
\scoH^{n+2}\bpar{P_n(\cbc), \pi_{n+1}(\cbc)}.
$$

The topological Postnikov invariant of a crossed complex $\cbc$ is the
image by this map of the algebraic Postnikov invariant $k_{n+1}\in H^2_{\cotrip[n]}
\bpar{P_n(\cbc), \tilde\pi^{(n)}_{n+1}(\cbc)}$.

For any space $X$ having the homotopy type of a crossed complex, we can obtain its
\pin s\ by simply calculating the topological \pin s\ of the fundamental
crossed complex of is singular complex  $\cbc = \Pi(X)$.

\appendix
\section{2-Torsors and Cotriple Cohomology} \label{apen}

Torsors were developed by Duskin as the appropriate algebraic structure
``representing'' general cotriple cohomology. The main references for this
subject are \cite{Duskin1975} and \cite{Glenn1982}. Our definitions differ
slightly from those found in those references in the sense that we put special
emphasis in the groupoid fibre of a torsor. This is closer to the way torsors
are defined and used, for example, in \cite{BuCaFa1998}.

\subsection{2-Torsors}

Every arrow in a groupoid establishes a group isomorphism between the
endomorphism group of its domain and that of its codomain.  A \emph{connected}
2-torsor with coefficients on a given abstract group $G$ is a connected groupoid
(meaning that it has one single connected component) together with a ``coherent"
system of isomorphisms from the different groups of endomorphisms to the
abstract group $G$. Thus, in order to specify a connected 2-torsor with
coefficients in $G$ we must give a connected groupoid \g\ (called ``\emph{the
fiber}" of the 2-torsor) and a natural system
$\alpha=\{\alpha_x\}_{x\in\obj(\g)}$ of group isomorphisms
$\alpha_x:\End_\g(x)\to G$.

This definition can be easily generalized to 2-torsors with a non-necessarily
connected fiber groupoid. As it turns out, the general definition can be
obtained as a particular case of the above, provided it is expressed in such a
way that it makes sense in more general categories.

Let \e\ be a Barr exact category. Associated to \e\ we have the category
$\Gpd(\e)$ of internal groupoid objects in \e\ and internal functors between
them. If $s,t:M\to O$ are the ``source" and ``target" structural morphisms of an
internal groupoid \g, we say that \g\ is connected if the coequalizer of $s$ and
$t$ is the terminal object, and we say that \g\ is totally disconnected if $s=t$.

For a given object $O$ in \e, $\tdGpd_O(\e)$ denotes the category whose objects
are totally disconnected internal groupoids in \e\ having $O$ as object of
objects, and whose arrows are internal functors whose component at the level of
objects is the identity of $O$. This category can be identified with the
category of internal group objects in the slice category $\e/O$.

If \g\ is an internal groupoid in \e, having $O$ as object of objects, the
category $\Gr(\e)^\g$ of internal \g-groups in \e\ is now defined in terms of
\g-actions. Then an internal \g-group consists of a totally disconnected
groupoid $\h\in\tdGpd_O$ together with a \g-action, that is a map in \e:
$$ M\times_O H\to H;\quad (f,h)\mapsto\;  \laction fh,
$$ where $M$ and $H$ denote the objects of arrows of \g\ and \h\ respectively,
and $M\times_O H$ is the pullback object of the diagram $M\xto{s}O
\xleftarrow{s=t}H$, such that it satisfies the usual axioms for an action of
groups. A morphism of internal \g-groups $\alpha:\h\to\h'$ is an
equivariant functor in $\tdGpd_O(\e)$.

One of the basic examples of internal \g-group is the \g-group of endomorphisms
of \g, $\End_\g$, defined by the totally disconnected groupoid $\End(\g)$ (the
equalizer of $s$ and $t$ with the groupoid structure given by restriction of
that in \g) and the action by conjugation in \g.

If we write $\mathbf{1}$ for the internal groupoid in \e, having the terminal
object as both object of arrows and object of objects, we can identify the
category $\Gr(\e)^\mathbf{1}$ with $\Gr(\e)$ and therefore the canonical $\g\to
\mathbf{1}$ induces a functor $\Gr(\e)\to\Gr(\e)^\g$ which allows us to regard
any internal group $G$ in $\e$ as a $\g$-group (the trivial action of $\g$ on
$G$).

Note that any internal functor $f:\g'\to\g$ defines, in a standard way, a
functor
$$ f^\ast:\Gr(\e)^\g\to\Gr(\e)^{\g'}.
$$

A \emph{connected 2-torsor} in \e\ consists then of a triple $(\g,G,\alpha)$
where $G$ is an internal group in \e, called the coefficients, \g\ is a
connected groupoid in \e, called the fiber, and $\alpha$, the cocycle, is an
isomorphism in the category $\Gr(\e)^\g$ from $\End_\g$ to $G$ (the trivial
action of \g\ on $G$). Note that to give the cocycle $\alpha$ is equivalent to
giving an arrow $\alpha:\End(\g)\to G$ that makes the following square a
pullback in \e,
$$
\xymatrix@C=2.5pc@R=1.75pc@!=1em {\End(\g)\ar[d]_{s\,=\, t}\ar[r]^-\alpha & G \ar[d]\\
\obj(\g)\ar[r] & \mathbf{1}}
$$ and satisfies:
\begin{itemize}
\item $\alpha(a b)=\alpha(a) \alpha(b)$ for all composable endomorphisms $a,b$
of \g, \item $\alpha(\laction fa)=\alpha(f a f^{-1})= \alpha(a)$, for all arrows
$f$ and all endomorphisms $a$ in \g, such that $s(f)=s(a)$.
\end{itemize}

The connected 2-torsors in \e\ whose group of coefficients is $G$ are the
objects of a category, denoted $\Tor^2(1,G)$, whose arrows form $(\g,G,\alpha)$
to $(\g',G,\alpha')$ are internal functors $f:\g\to\g'$
compatible with the cocycles $\alpha,\alpha'$ in the sense that
$$
\alpha = f^\ast(\alpha')= \alpha'\circ f.
$$

If $T$ is an object in a Barr exact category \e, then the slice category $\e/T$
is again a Barr exact category and a \emph{2-torsor above} $T$ in \e\ is defined
as a connected 2-torsor in $\e/T$. In this definition it is understood that the
coefficients are taken in an internal group object  in $\e / T$. If $T$ is an
object and $G$ an internal group in \e, then the canonical projection $G\times
T\to T$ gives an object of $\e/T$ which has a canonical structure of internal
group object in $\e/T$. In this situation, a $(G,2)$-torsor above $T$ in \e\ is
defined as a connected 2-torsor in $\e/T$ with coefficients in $G\times T$, and
the category of such torsors is denoted $\Tor^2(T,G)$. By $\Tor^2[T,G]$ we
denote the set of connected components of $\Tor^2(T,G)$.

Let us suppose now that we have a functor $U:\e\to\s$ from a Barr exact category
to a category with finite limits. Let \g\ be an internal groupoid in \e\ with
source and target maps $s,t:M\to O$. If $q:O\to T$ is the coequalizer of $s$ and
$t$ and $O\times_T O$ is the pullback of $q$ with itself, there is an induced
map $(s,t):M\to O\times_T O$. We say that the groupoid \g\ is \emph{$U$-split}
if the maps $U(q)$ and $U(s,t)$ split in \s.

A $U$-\emph{split} $(G,2)$-torsor above $T$ is a $(G,2)$-torsor above $T$ such
that its fiber groupoid is $U$-split. The full subcategory of $\Tor^2(T,G)$
determined by those $(G,2)$-torsors which are $U$-split is denoted
$\Tor_U^2(T,G)$. Correspondingly, the set of connected components of
$\Tor_U^2(T,G)$ is denoted $\Tor_U^2[T,G]$.

\guardar{If $(\g,\alpha)$ is a $(G,2)$-torsor above $T$ and $f:O'\to
O=\obj(\g)$ is a morphism in \e, then there is a $(G,2)$-torsor
$(\g',\alpha')$ above $T$ with a cartesian morphism of torsors
$f':(\g',\alpha')\to(\g,\alpha)$ whose image by
$\obj:\Gpd(\e)\to\e$ is $f$}

\begin{proposition}
\label{condition U-split} $\Tor^2(T,G)$ is a category fibred over $\e$ via the
composite functor
$$
\everyentry={\vphantom{\Big(}} \xymatrix@C=1.5pc@R=1.75pc
{\Tor^2(T,G)\ar[r]^-{\fib} & \Gpd(\e)\ar[r]^-{\obj} & \e.}
$$ Furthermore, if $(\g,\alpha)\in \Tor^2(T,G)$ is $U$-split for some left exact
functor $U:\e\to\s$, then $(\g',\alpha')$ is $U$-split if and only
if the projection $q':O'\to T$ corresponding to
$(\g',\alpha')$ is $U$-split.
\end{proposition}

\begin{proof} For the first part we have to prove that if $(\g,\alpha)$ is a
$(G,2)$-torsor above $T$ and $f:O'\to O=\obj(\g)$ is a morphism in \e,
then there is a $(G,2)$-torsor $(\g',\alpha')$ above $T$ and an
$(\obj\circ\fib)$-cartesian morphism of torsors
$f':(\g',\alpha')\to(\g,\alpha)$ above $f$.

The idea for the construction of $\g'$ with object of objects $O'$
is that its arrows from $x\in O'$ to $y\in O'$ are the arrows
$f(x)\to f(y)$ in \g\ and that the identity of $x\in O'$ is the identity
of $f(x)$. The composition in $\g'$ will then be clearly induced by that
of \g. Thus, the object of arrows of $\g'$, together with its source and
target maps can internally be defined by the following pullback
$$
\xymatrix@C=1.6pc@R=1.75pc@!=1em {M'\ar[d]_u \ar[rr]^-{(s',t')} & & O'\times_T O' \ar[d]^{f\times
f}
\\ M\ar[rr]_-{(s,t)} && O\times_T O'}
$$ where $O\times_TO$ is the pullback of $q$ with itself and
$O'\times_TO'$ that of $q'=qf$ with itself. This construction
produces a functor $f':\g'\to\g$ whose component on arrows is $u$
(and it is $f$ on objects). Note that $\End_{\g'}=
f^{\prime\ast}(\End_\g)$. The cocycle map $\alpha'$ is defined as the
image of $\alpha$ by the induced functor
$f^{\prime\ast}:\Gr(\e)^\g\to\Gr(\e)^{\g'}$, that is
$\alpha'=f^{\prime\ast}\alpha$.

Let us now assume that $(\g,\alpha)$ is $U$-split. For
$(\g',\alpha')$ to be $U$-split it is necessary that $q'$ be
$U$-split. From the exactness of $U$ and a splitting of $U(s,t)$ it follows that
$(s',t')$ is $U$-split, therefore that $q'$ be $U$-split is
also sufficient for $(\g',\alpha')$ to be $U$-split.

It only remains to prove that $f'$ is cartesian. Let $g:\g''\to\g$ be a
morphism of internal groupoid in \e\ such that $g_0=\obj(g):O''\to O$ factors
through $f$ as $g_0=fh$. Then it is clear which is the only way to define a factorization
$g=f' h'$ such that $\obj(h')=h$. This condition determines $h'$ on
objects and $g$ determines it on arrows.
\end{proof}

By a reasoning similar to the one given in \cite{Glenn1982} to prove Theorem 5.7.5,
it is easy to prove that from any diagram in $\Tor_U^2(T,G)$ of the form
$(\g,\alpha) \to (\tilde\g,\tilde\alpha)\lto(\g',\alpha')$ one can obtain another of
the form $ (\g,\alpha) \lto (\g'', \alpha'') \to (\g',\alpha')$ (see
\cite{Garcia2003}, Lema 4.2.6 for the details). As a consequence we have the
following useful necessary (and, obviously, also sufficient) condition satisfied by
torsors in the same connected component of $\Tor_U^2(T,G)$.

\begin{proposition}
\label{adapted from Glenn} If $(\g,\alpha)$ and $(\g',\alpha')$ are
$U$-split 2-torsors which are in the same connected component of
$\Tor_U^2(T,G)$, then there is a torsor
$(\g'',\alpha'')$ and maps
$$ (\g,\alpha) \lto
(\g'',\alpha'')\to(\g',\alpha')
$$ in $\Tor_U^2(T,G)$.
\end{proposition}


\subsection{Cotriple Cohomology}

Let \e\ be a tripleable category over a category \s\ with cotriple \bbg. For any
object $T\in\e$ and any abelian group object $A$ in $\e/T$, the cotriple
cohomology groups $H^n_{\bbg}(T,A)$ can be represented in terms of homotopy
classed of simplicial maps from the cotriple simplicial resolution of $T$ to the
Eilenberg-Mac Lane complex $K(A,n)$. On the other hand, Duskin's interpretation
theorem for cotriple cohomology \cite{Duskin1975} provides an interpretation of
the elements of $H^n_{\bbg}(T,A)$ in terms of $U$-split torsors, where
$U:\e\to\s$ is the monadic functor associated to the cotriple \bbg. In the
particular case $n=2$, Duskin's theorem implies the following

\begin{theorem}
\label{Duskin's interpretation} In the above conditions, for any object $T\in\e$
and any abelian group object $A$ in $\e/T$, there is a natural bijection
$$ H^2_{\bbg}(T,A)\cong \Tor^2_U[T,A].
$$
\end{theorem}

In the presence of a cotriple, Proposition \ref{condition U-split} has the
following consequence:

\begin{proposition}
\label{cotriple torsor in same component} In the hypothesis of Theorem
\ref{Duskin's interpretation}, if $(\g,\alpha)$ is a $U$-split $(A,2)$-torsor
above $T$ with object of objects $O$, there is a $U$-split $(A,2)$-torsor above
$T$, $(\g',\alpha')$, whose object of objects is $\bbg(T)$ and whose
projection is the counit $\varepsilon_T:\bbg(T)\to T$. Furthermore,
$(\g',\alpha')$ is connected to $(\g,\alpha)$ by a morphism
$(\g',\alpha')\to(\g,\alpha)$ in $\Tor^2_U(T,A)$.
\end{proposition}

\begin{proof} Let $s:U(T)\to U(O)$ be a section of the image $U(q)$ of the
projection $q:O\to T$ of $(\g,\alpha)$. Use Proposition \ref{condition U-split}
with $O'=\bbg(T)$ and $f$ equal to the composite
$\bbg(T)\xrightarrow{F(s)}\bbg(O)\xrightarrow{\varepsilon_0} O$, where $F$ is
the left adjoint to $U$ and $\varepsilon$ is the counit of \bbg\ ($=FU$). Then
we obtain an $(A,2)$-torsor above $T$, $(\g',\alpha')$, whose
projection is the composite
$$ qf=q\varepsilon:OF(s)= \varepsilon_T FU(q)F(s)=\varepsilon_T,
$$ and a map $f':(\g',\alpha')\to(\g,\alpha)$ which is given by $f$ at the level
of objects. Since $\varepsilon_T$ is a $U$-split map (with $U$-section given by
$\eta_{U(O)}$ ), it follows that $(\g',\alpha')$ is $U$-split.
\end{proof}

Let now $\utor_U^2(T,A)$ denote the full subcategory of $\Tor_U^2(T,A)$
determined by those torsors whose object of objects is $\bbg(T)$ and whose
projection is the counit $\varepsilon_T$. Then Propositions \ref{cotriple torsor
in same component} and \ref{adapted from Glenn} imply the following:

\begin{proposition}
\label{same connected components} In the hypothesis of Theorem \ref{Duskin's
interpretation}, let $F:\be\to\Tor_U^2(T,A)$ be a full and faithful functor such
that the inclusion $\utor_U^2(T,A)\hto\Tor_U^2(T,A)$ factors through $F$. Then,
$F$ establishes a bijection between the set $[\be]$ of connected components of
$\be$ and  $\Tor_U^2[T,A]$. Hence there is a natural bijection
$$ H^2_\bbg(T,A)\cong [\be].
$$
\end{proposition}

\begin{proof} Let $A,B\in\be$ such that $F(A)$ and $F(B)$ are in the same
connected component of $\Tor_U^2(T,A)$. We just need to show that $A$ and $B$ are in the same
connected component in $\be$. By Proposition \ref{cotriple torsor in same component} there are
2-torsors $X, Y\in \utor_U^2(T,A)$ and morphisms $h:X\to F(A)$ and $k:Y\to F(B)$, such that $X$
and $Y$ are in the same connected component of $\Tor_U^2(T,A)$. Hence, by Proposition
\ref{adapted from Glenn} we get a diagram $\xymatrix@1@C=1pc{X&Z\labelmargin{1pt} \ar[l]_(.45)f
\labelmargin{1.5pt}\ar[r]^g &Y}$ in $\Tor_U^2(T,A)$ where, by Proposition \ref{cotriple torsor in
same component} we can suppose that $Z$ is in $\utor_U^2(T,A)$. By the hypothesis that the
inclusion of $\utor_U^2(T,A)$ factors through $F$, we get a diagram
\begin{equation}
\label{first diagram}A' \lto C \to B'
\end{equation} in $\be$ such that $F(A') = X$ and $F(B') = Y$. Thus, the maps
$h, k$ and the fullness of $F$ allow us to extend diagram \eqref{first diagram}
to a diagram
$$ A \lto A' \lto C \to B' \to B
$$ proving that $A$ and $B$ are in the same connected component.
\end{proof}

\end{document}